\documentclass[11pt]{article}
\usepackage{setspace}
\usepackage{amsmath} %provare a modificare per toglierli
\usepackage{amssymb}
\usepackage{latexsym}
\usepackage{xcolor}
\usepackage{fancyhdr}
\allowdisplaybreaks 

 \usepackage{comment}

\def\XXint#1#2#3{{\setbox0=\hbox{$#1{#2#3}{\int}$} 
\vcenter{\hbox{$#2#3$}}\kern-.5\wd0}}   
 \numberwithin{equation}{section}
\newtheorem{thm}[equation]{Theorem}
\newtheorem{prop}[equation]{Proposition}
\newtheorem{defn}[equation]{Definition}
\newtheorem{rem}[equation]{Remark}
\newtheorem{lem}[equation]{Lemma}

\title{
Continuity  of the double layer potential
 \\
of a second order elliptic differential operator 
\\
 in Schauder spaces on the boundary} 
 
\author{  
Massimo Lanza de Cristoforis
\\
Dipartimento di Matematica `Tullio Levi-Civita', 
\\
Universit\`a degli Studi di Padova, 
\\
Via Trieste 63, Padova 35121, 
Italy. 
\\
E-mail: mldc@math.unipd.it   }

\date{\ }

\begin{document}
 \maketitle

%\noindent
%{\bf Running title:} \\

\noindent
{\bf Abstract:}  We prove the validity of a regularizing property on the boundary of the  double layer potential associated to the fundamental solution of a {\em nonhomogeneous} second order elliptic differential operator with constant coefficients in Schauder spaces of exponent greater or equal to two that sharpens classical results of N.M.~G\"{u}nter, S.~Mikhlin, V.D.~Kupradze, T.G.~Gegelia, M.O.~Basheleishvili and T.V.~Bur\-chuladze, U.~Heinemann and extends the work of  
A.~Kirsch who has considered the case of the Helmholtz operator. 
 
 \vspace{\baselineskip}

\par


\noindent
{{\bf ORCID iD:}} 0000 0001 6886 4647
\section{Introduction} In this paper, we consider the double layer potential associated to the fundamental solution of a second order differential operator with constant coefficients. Unless otherwise specified, we assume that
\[
n\in {\mathbb{N}}\setminus\{0,1\}\,,
\]
where ${\mathbb{N}}$ denotes the set of natural numbers including $0$. Let $\alpha\in]0,1]$, $m\in {\mathbb{N}}\setminus\{0\}$. Let $\Omega$ be a bounded open subset of ${\mathbb{R}}^{n}$ of class $C^{m,\alpha}$. For the notation and standard properties of the H\"{o}lder and of the Schauder spaces $C^{m,\alpha}$ we refer  to \cite[\S 2]{DoLa17}, \cite[\S 2.6,  2.11]{DaLaMu21}. Let $\nu \equiv (\nu_{l})_{l=1,\dots,n}$ denote the external unit normal to $\partial\Omega$. Let $N_{2}$ denote the number of multi-indexes $\gamma\in {\mathbb{N}}^{n}$ with $|\gamma|\leq 2$. For each 
\begin{equation}
\label{introd0}
{\mathbf{a}}\equiv (a_{\gamma})_{|\gamma|\leq 2}\in {\mathbb{C}}^{N_{2}}\,, 
\end{equation}
we set 
\[
a^{(2)}\equiv (a_{lj} )_{l,j=1,\dots,n}\qquad
a^{(1)}\equiv (a_{j})_{j=1,\dots,n}\qquad
a\equiv a_{0}\,.
\]
with $a_{lj} \equiv 2^{-1}a_{e_{l}+e_{j}}$ for $j\neq l$, $a_{jj} \equiv
 a_{e_{j}+e_{j}}$,
and $a_{j}\equiv a_{e_{j}}$, where $\{e_{j}:\,j=1,\dots,n\}$  is the canonical basis of ${\mathbb{R}}^{n}$. We note that the matrix $a^{(2)}$ is symmetric. 
Then we assume that 
  ${\mathbf{a}}\in  {\mathbb{C}}^{N_{2}}$ satisfies the following ellipticity assumption
\begin{equation}
\label{ellip}
\inf_{
\xi\in {\mathbb{R}}^{n}, |\xi|=1
}{\mathrm{Re}}\,\left\{
 \sum_{|\gamma|=2}a_{\gamma}\xi^{\gamma}\right\} >0\,,
\end{equation}
and we consider  the case in which
\begin{equation}
\label{symr}
a_{lj} \in {\mathbb{R}}\qquad\forall  l,j=1,\dots,n\,.
\end{equation}
Then we introduce the operators
\begin{eqnarray*}
%\label{introd1}
P[{\mathbf{a}},D]u&\equiv&\sum_{l,j=1}^{n}\partial_{x_{l}}(a_{lj}\partial_{x_{j}}u)
+
\sum_{l=1}^{n}a_{l}\partial_{x_{l}}u+au\,,
\\
%\label{introd2}
B_{\Omega}^{*}v&\equiv&\sum_{l,j=1}^{n} \overline{a}_{jl}\nu_{l}\partial_{x_{j}}v
-\sum_{l=1}^{n}\nu_{l}\overline{a}_{l}v\,,
\end{eqnarray*}
for all $u,v\in C^{2}(\overline{\Omega})$, and a fundamental solution $S_{{\mathbf{a}} }$ of $P[{\mathbf{a}},D]$, and the  boundary integral operator corresponding to the  double layer potential 
\begin{eqnarray}
\label{introd3}
\lefteqn{
 W_\Omega[{\mathbf{a}},S_{{\mathbf{a}}}   ,\mu](x) \equiv 
 \int_{\partial\Omega}\mu (y)\overline{B^{*}_{\Omega,y}}\left(S_{{\mathbf{a}}}(x-y)\right)
\,d\sigma_{y}
}
\\  \nonumber
&&
\qquad\qquad\qquad
=- \int_{\partial\Omega}\mu(y)\sum_{l,j=1}^{n} a_{jl}\nu_{l}(y)\frac{\partial S_{ {\mathbf{a}} } }{\partial x_{j}}(x-y)\,d\sigma_{y}
\\  \nonumber
&&
\qquad\quad\qquad\quad
-\int_{\partial\Omega}\mu(y)\sum_{l=1}^{n}\nu_{l}(y)a_{l}
S_{ {\mathbf{a}} }(x-y)\,d\sigma_{y}  
\end{eqnarray}
for all $x\in \partial\Omega$,
where the density or moment $\mu$ is a function  from $\partial\Omega$ to ${\mathbb{C}}$
and $d\sigma_{y}$ is the ordinary $(n-1)$-dimensional measure.
 Here the subscript $y$ of $\overline{B^{*}_{\Omega,y}}$ means that we are taking $y$ as variable of the differential operator $\overline{B^{*}_{\Omega,y}}$. 
The role of the double layer potential in the solution of boundary value problems for the operator $P[{\mathbf{a}},D]$ is well known (cf.~\textit{e.g.}, 
G\"{u}nter~\cite{Gu67}, Kupradze,  Gegelia,  Basheleishvili and 
 Burchuladze~\cite{KuGeBaBu79}, Mikhlin \cite{Mik70}, 
Mikhlin and  Pr\"{o}ssdorf \cite{MikPr86}, Buchukuri,  Chkadua,  Duduchava, and  Natroshvili \cite{BuChDuNa12}.) 
 
We now briefly summarize some known results in the classical case of the boundary behaviour of the  double layer potential in Schauder spaces with $m\geq 2$. Instead for the regularity properties of the double layer potential in Schauder spaces with $m\geq 2$ outside of the boundary we refer to G\"{u}nter~\cite{Gu67}, Kupradze,  Gegelia,  Basheleishvili and 
 Burchuladze~\cite{KuGeBaBu79}, Mikhlin \cite{Mik70}, 
Mikhlin and  Pr\"{o}ssdorf \cite{MikPr86}, Miranda \cite{Mi65},   \cite{Mi70},  
Wiegner~\cite{Wi93}, Dalla Riva \cite{Da13}, Dalla Riva, Morais and Musolino \cite{DaMoMu13}, Mitrea, Mitrea and Verdera \cite{MitMitVe16} and references therein.

In case  $n=3$, $m\geq 2$, $\alpha\in]0,1]$ and $\Omega$ is of class $C^{m,\alpha}$ and if  $P[{\mathbf{a}},D]$ is the  Laplace operator,  G\"{u}nter~\cite[Appendix, \S\ IV, Thm.~3]{Gu67} has 
proved that   $W[\partial\Omega ,{\mathbf{a}},S_{{\mathbf{a}}},\cdot]$ is bounded from $C^{m-2,\alpha}(\partial\Omega)$ to $C^{m-1,\alpha'}(\partial\Omega)$ for $\alpha'\in]0,\alpha[$.

In case $n\geq 2$, $m\geq 2$, $\alpha\in]0,1]$,  O.~Chkadua \cite{Chk23} has pointed out  that one could exploit Kupradze,  Gegelia,  Basheleishvili and 
 Burchuladze~\cite[Chap.~IV, Sect.~2, Thm 2.9, Chap. IV, Sect.~3, Theorems 3.26 and  3.28]{KuGeBaBu79} and  prove that if $\Omega$ is of class $C^{m,\alpha}$, then $W[\partial\Omega ,{\mathbf{a}},S_{{\mathbf{a}}},\cdot]$ is bounded from $C^{m-1,\alpha'}(\partial\Omega)$ to $C^{m,\alpha'}(\partial\Omega)$ for $\alpha'\in]0,\alpha[$.\par

In case $n=3$  and $\Omega$ is of class $C^{2}$, $\alpha\in]0,1[$   and if  $P[{\mathbf{a}},D]$ is the  Helmholtz operator, Colton and Kress~\cite{CoKr83} have 
developed previous work of G\"{u}nter~\cite{Gu67} and Mikhlin~\cite{Mik70}  and proved that the operator $W[\partial\Omega ,{\mathbf{a}},S_{{\mathbf{a}}},\cdot]$
is bounded from $C^{0,\alpha}(\partial\Omega)$ to $C^{1,\alpha}(\partial\Omega)$.

In case $n\geq 2$, $\alpha\in]0,1[$  and $\Omega$ is of class $C^{2}$   and if  $P[{\mathbf{a}},D]$ is the  Laplace operator, 
Hsiao and Wendland \cite[Remark 1.2.1]{HsWe08} 
deduce  that the operator $W[\partial\Omega ,{\mathbf{a}},S_{{\mathbf{a}}},\cdot]$
is bounded from $C^{0,\alpha}(\partial\Omega)$ to $C^{1,\alpha}(\partial\Omega)$ by the work of
Mikhlin and Pr\"{o}ssdorf \cite{MikPr86}.

In case  $n=3$, $m\geq 2$ and $\Omega$ is of class $C^{m,\alpha}$  and if  $P[{\mathbf{a}},D]$ is the  Helmholtz operator, Kirsch~\cite[Thm.~3.3 (a)]{Ki89} has 
developed previous work of  G\"{u}nter~\cite{Gu67}, Mikhlin~\cite{Mik70} and Colton and Kress~\cite{CoKr83} 
and has proved that  the operator  $W_\Omega[{\mathbf{a}},S_{{\mathbf{a}}}   ,\cdot]$ is  bounded from $C^{m-1,\alpha}(\partial\Omega)$ to $C^{m,\alpha}(\partial\Omega)$. 
 
 von Wahl~\cite{vo90} has considered the case of Sobolev spaces and has proved that if $\Omega$ is of class $C^{\infty}$ and if 
 $S_{{\mathbf{a}}}$ is the fundamental solution of the Laplace operator,  then  the double layer improves the regularity of one unit on the boundary. Then Heinemann~\cite{He92} has developed the ideas of  von Wahl in the frame of Schauder spaces and has proved that
 if $\Omega$ is of class $C^{m+5}$ and if 
 $S_{{\mathbf{a}}}$ is the fundamental solution of the Laplace operator,  then  the double layer improves the regularity of one unit on the boundary, \textit{i.e.}, 
 $W_\Omega[{\mathbf{a}},S_{{\mathbf{a}}}   ,\cdot]$ is linear and continuous from $C^{m,\alpha}(\partial\Omega)$ to $C^{m+1,\alpha}(\partial\Omega)$.

Maz'ya and  Shaposhnikova~\cite{MaSh05}
 have proved that $W_\Omega[{\mathbf{a}},S_{{\mathbf{a}}}   ,\cdot]$ is continuous in fractional Sobolev spaces under sharp   regularity assumptions on the boundary and if  $P[{\mathbf{a}},D]$ is the  Laplace operator. 
 
 Dondi and the author  \cite{DoLa17}   have proved that  if $m\geq 2$ and $\Omega$ is of class $C^{m,\alpha}$ with $\alpha\in]0,1[$, then the double layer potential $W_\Omega[{\mathbf{a}},S_{{\mathbf{a}}}   ,\cdot]$ associated to the fundamental solution of a {\em nonhomogeneous} second order elliptic differential operator with constant coefficients
 is  bounded from $C^{m,\beta}(\partial\Omega)$ to $C^{m,\alpha}(\partial\Omega)$
 for all $\beta\in]0,\alpha]$. 
 
 For corresponding results for the fundamental solution of the heat equation, we refer to the  author and Luzzini   \cite{LaLu17}, \cite{LaLu18}  and references therein.
 
 In this paper we plan to prove that  if $m\geq 2$ and $\Omega$ is of class $C^{m,\alpha}$ with $\alpha\in]0,1]$, then the double layer potential $W_\Omega[{\mathbf{a}},S_{{\mathbf{a}}}   ,\cdot]$ associated to the fundamental solution of a {\em nonhomogeneous} second order elliptic differential operator with constant coefficients
 is  bounded from $C^{m-1,\alpha}(\partial\Omega)$ to $C^{m,\alpha}(\partial\Omega)$ in case $\alpha<1$ and to the generalized Schauder space
$C^{m,\omega_{1}(\cdot)}(\partial\Omega)$ of functions with $m$-th order derivatives which satisfy a generalized $\omega_{1}(\cdot)$-H\"{o}lder condition with $\omega_{1}(\cdot)$ as in (\ref{omth}) and thus with 
\[
\omega_{1}(r)\sim r|\ln r| 
\qquad{\mathrm{as}}\ r\to 0,
\]
 in case $\alpha=1$. For the classical definition of the generalized H\"{o}lder or Schauder spaces on the boundary, we refer the reader to the author and Dondi  
  \cite[\S 2]{DoLa17}   and to Dalla Riva, the author  and Musolino \cite[\S 2.6,   2.20]{DaLaMu21}.\par

  Hence we sharpen the work of the above mentioned authors in the sense that if 
  $\Omega$ is of class $C^{m,\alpha}$ with $m\geq 2$, then the class of regularity of  the target space of $W_\Omega[{\mathbf{a}},S_{{\mathbf{a}}}   ,\cdot]$ is precisely $C^{m,\alpha}$  if  $\alpha<1$
  and is the generalized Schauder space $C^{m,\omega_1(\cdot)}$ if $\alpha=1$.
  
Moreover, we extend 
  the above mentioned result of  Kirsch~\cite{Ki89}  in the sense that Kirsch~\cite{Ki89}
   has  considered the Helmholtz operator in case $n=3$, $\alpha<1$ and we have considered a general fundamental solution $S_{{\mathbf{a}}}$ with ${\mathbf{a}}$ as in (\ref{introd0}), (\ref{ellip}), (\ref{symr}), $\alpha\leq 1$   and $n\geq 2$.
 
  \section{Notation}\label{sec:notation}
  Let $M_n({\mathbb{R}})$ denote the set of $n\times n$ matrices with real entries. $|A|$ denotes the operator norm of a matrix $A$, 
       $A^{t}$ denotes the transpose matrix of $A$. Let $O_{n}({\mathbb{R}})$ denote the set of $n\times n$ orthogonal matrices with real entries.  We set
\begin{equation}\label{eq:balls}
{\mathbb{B}}_n(\xi,r)\equiv \left\{\eta\in {\mathbb{R}}^n:\, |\xi-\eta|<r\right\}
\qquad \forall
(\xi,r)\in {\mathbb{R}}^n\times ]0,+\infty[
\,.
\end{equation}
 If ${\mathbb{D}}$ is a subset of $ {\mathbb{R}}^n$, 
then we set
\[
B({\mathbb{D}})\equiv\left\{
f\in {\mathbb{C}}^{\mathbb{D}}:\,f\ \text{is\ bounded}
\right\}
\,,\quad
\|f\|_{B({\mathbb{D}})}\equiv\sup_{\mathbb{D}}|f|\qquad\forall f\in B({\mathbb{D}})\,.
\]
Then $C^0({\mathbb{D}})$ denotes the set of continuous functions from ${\mathbb{D}}$ to ${\mathbb{C}}$ and we introduce the subspace
$
C^0_b({\mathbb{D}})\equiv C^0({\mathbb{D}})\cap B({\mathbb{D}})
$
of $B({\mathbb{D}})$.  Let $\omega$ be a function from $[0,+\infty[$ to itself such that
\begin{eqnarray}
\nonumber
&&\qquad\qquad\omega(0)=0,\qquad \omega(r)>0\qquad\forall r\in]0,+\infty[\,,
\\
\label{om}
&&\qquad\qquad\omega\ {\text{is\   increasing,}}\ \lim_{r\to 0^{+}}\omega(r)=0\,,
\\
\nonumber
&&\qquad\qquad{\text{and}}\ \sup_{(a,t)\in[1,+\infty[\times]0,+\infty[}
\frac{\omega(at)}{a\omega(t)}<+\infty\,.
\end{eqnarray}
Here `$\omega$ is increasing' means that 
$\omega(r_1)\leq \omega(r_2)$ whenever $r_1$, $r_2\in [0,+\infty[$ and $r_1<r_2$.
If $f$ is a function from a subset ${\mathbb{D}}$ of ${\mathbb{R}}^n$   to ${\mathbb{C}}$,  then we denote by   $|f:{\mathbb{D}}|_{\omega (\cdot)}$  the $\omega(\cdot)$-H\"older constant  of $f$, which is delivered by the formula   
\[
%\label{om1a}
|f:{\mathbb{D}}|_{\omega (\cdot)
}
\equiv
\sup\left\{
\frac{|f( x )-f( y)|}{\omega(| x- y|)
}: x, y\in {\mathbb{D}} ,  x\neq
 y\right\}\,.
\]        
If $|f:{\mathbb{D}}|_{\omega(\cdot)}<\infty$, we say that $f$ is $\omega(\cdot)$-H\"{o}lder continuous. Sometimes, we simply write $|f|_{\omega(\cdot)}$  
instead of $|f:{\mathbb{D}}|_{\omega(\cdot)}$. The
subset of $C^{0}({\mathbb{D}} ) $  whose
functions  are
$\omega(\cdot)$-H\"{o}lder continuous    is denoted  by  $C^{0,\omega(\cdot)} ({\mathbb{D}})$
and $|f:{\mathbb{D}}|_{\omega(\cdot)}$ is a semi-norm on $C^{0,\omega(\cdot)} ({\mathbb{D}})$.  
Then we consider the space  $C^{0,\omega(\cdot)}_{b}({\mathbb{D}} ) \equiv C^{0,\omega(\cdot)} ({\mathbb{D}} )\cap B({\mathbb{D}} ) $ with the norm \[
%\label{om1b}
\|f\|_{ C^{0,\omega(\cdot)}_{b}({\mathbb{D}} ) }\equiv \sup_{x\in {\mathbb{D}} }|f(x)|+|f|_{\omega(\cdot)}\qquad\forall f\in C^{0,\omega(\cdot)}_{b}({\mathbb{D}} )\,.
\] 
\begin{rem}
\label{rem:om4}
Let $\omega$ be as in (\ref{om}). 
Let ${\mathbb{D}}$ be a   subset of ${\mathbb{R}}^{n}$. Let $f$ be a bounded function from $ {\mathbb{D}}$ to ${\mathbb{C}}$, $a\in]0,+\infty[$.  Then,
\[
\label{rem:om5}
\sup_{x,y\in {\mathbb{D}},\ |x-y|\geq a}\frac{|f(x)-f(y)|}{\omega(|x-y|)}
\leq \frac{2}{\omega(a)} \sup_{{\mathbb{D}}}|f|\,.
\]
\end{rem}
In the case in which $\omega(\cdot)$ is the function 
$r^{\alpha}$ for some fixed $\alpha\in]0,1]$, a so-called H\"{o}lder exponent, we simply write $|\cdot:{\mathbb{D}}|_{\alpha}$ instead of
$|\cdot:{\mathbb{D}}|_{r^{\alpha}}$, $C^{0,\alpha} ({\mathbb{D}})$ instead of $C^{0,r^{\alpha}} ({\mathbb{D}})$, $C^{0,\alpha}_{b}({\mathbb{D}})$ instead of $C^{0,r^{\alpha}}_{b} ({\mathbb{D}})$, and we say that $f$ is $\alpha$-H\"{o}lder continuous provided that 
$|f:{\mathbb{D}}|_{\alpha}<\infty$. For the standard properties of the spaces of H\"{o}lder or Lipschitz continuous functions, we refer    to \cite[\S 2]{DoLa17}, \cite[\S 2.6]{DaLaMu21}.  Let $\Omega$ be an open subset of ${\mathbb{R}}^n$.  Let $s\in {\mathbb{N}}\setminus\{0\}$, $f\in \left(C^{1}(\Omega)\right)^{s} $. Then   $Df$ denotes the Jacobian matrix of $f$.

\section{Special classes of potential type kernels in ${\mathbb{R}}^n$}\label{sec:kecla}
 In this section we collect some basic properties of the classes of kernels that we need. For the proofs, we refer to \cite[\S 3]{La22b}. If $X$ and $Y$ are subsets of ${\mathbb{R}}^n$, then we denote by ${\mathbb{D}}_{X\times Y}$ the diagonal  of $X\times Y$, i.e., we set
 \begin{equation}\label{diagonal}
{\mathbb{D}}_{X\times Y}\equiv\left\{
(x,y)\in X\times Y:\,x=y
\right\} 
\end{equation}
and if $X=Y$, then    we denote by ${\mathbb{D}}_{X}$ the diagonal  of $X\times X$, i.e., we set
\[
{\mathbb{D}}_X\equiv {\mathbb{D}}_{X\times X}\,.
\]
An off-diagonal function  in $X\times Y$ is a function from $(X\times Y)\setminus {\mathbb{D}}_{X\times Y}$ to ${\mathbb{C}}$. We plan to consider the well known class of potential type  off-diagonal kernels as in the following definition.
\begin{defn}
 Let $X$ and $Y$ be subsets of ${\mathbb{R}}^n$. Let $s\in {\mathbb{R}}$. We denote by ${\mathcal{K}}_{s,X\times Y}$ (or more simply by ${\mathcal{K}}_s$), the set of continuous functions $K$ from $(X\times Y)\setminus {\mathbb{D}}_{ X\times Y }$ to ${\mathbb{C}}$ such that
\[
 \|K\|_{ {\mathcal{K}}_{s,X\times Y} }\equiv \sup_{(x,y)\in  (X\times Y)\setminus {\mathbb{D}}_{ X\times Y }  }|K(x,y)|\,|x-y|^s<+\infty\,.
\]
The elements of $ {\mathcal{K}}_{s,X\times Y}$ are said to be kernels of potential type $s$ in $X\times Y$. 
\end{defn}
We plan to consider specific classes of `potential type' kernels that are suitable to prove continuity theorems for integral operators in H\"{o}lder spaces
as in the following definition, which is a generalisation of related classes as in Gegelia,  Basheleishvili and  Burchuladze~\cite{KuGeBaBu79} (see also Dondi and the author  \cite{DoLa17},  where such classes have been introduced in a form that generalizes those of Giraud \cite{Gi34}, Gegelia \cite{Ge67} and Gegelia,  Basheleishvili and  Burchuladze~\cite[Chap.~IV]{KuGeBaBu79}).
\begin{defn}\label{defn:ksss}
 Let $X$, $Y\subseteq {\mathbb{R}}^n$. Let $s_1$, $s_2$, $s_3\in {\mathbb{R}}$. We denote by ${\mathcal{K}}_{s_1, s_2, s_3} (X\times Y)$ the set of continuous functions $K$ from $(X\times Y)\setminus {\mathbb{D}}_{X\times Y}$ to ${\mathbb{C}}$ such that
 \begin{eqnarray*}
\lefteqn{
\|K\|_{  {\mathcal{K}}_{ s_1, s_2, s_3  }(X\times Y)  }
\equiv
\sup\biggl\{\biggr.
|x-y|^{ s_{1} }|K(x,y)|:\,(x,y)\in X\times Y, x\neq y
\biggl.\biggr\}
}
\\ \nonumber
&&\qquad\qquad\qquad
+\sup\biggl\{\biggr.
\frac{|x'-y|^{s_{2}}}{|x'-x''|^{s_{3}}}
|  K(x',y)- K(x'',y)  |:\,
\\ \nonumber
&&\qquad\qquad\qquad 
x',x''\in X, x'\neq x'', y\in Y\setminus{\mathbb{B}}_{n}(x',2|x'-x''|)
\biggl.\biggr\}<+\infty\,.
\end{eqnarray*}
\end{defn}
 One can easily verify that $({\mathcal{K}}_{ s_{1},s_{2},s_{3}   }(X\times Y),\|\cdot\|_{  {\mathcal{K}}_{s_{1},s_{2},s_{3}   }(X\times Y)  })$ is a normed space.  By our definition, if $s_1$, $s_2$, $s_3\in {\mathbb{R}}$, we have
\[
{\mathcal{K}}_{s_{1},s_{2},s_{3}   }(X\times Y) \subseteq {\mathcal{K}}_{s_{1}, X\times Y} 
\]
and
\[
\|K\|_{{\mathcal{K}}_{s_{1}, X\times Y} }\leq \|K\|_{ {\mathcal{K}}_{s_{1},s_{2},s_{3}   }(X\times Y) }
\qquad\forall K\in {\mathcal{K}}_{s_{1},s_{2},s_{3}   }(X\times Y) \,.
\]
We note that if we choose $s_2=s_1+s_3$ we have a so-called class of standard kernels.   Then we have the following elementary known embedding lemma (cf.~\textit{e.g.}, \cite[Lem.~3.1]{La22b}).
\begin{lem}\label{lem:kelem}
 Let $X$, $Y\subseteq {\mathbb{R}}^n$. Let $s_1$, $s_2$, $s_3\in {\mathbb{R}}$. If $a\in]0,+\infty[$, then $ {\mathcal{K}}_{s_{1},s_{2},s_{3}   }(X\times Y)$ is continuously embedded into  ${\mathcal{K}}_{s_{1},s_{2}-a,s_{3}-a   }(X\times Y)$.
\end{lem}
Next we state the following two product rule statements (cf.~\cite[Thm.~3.1, Prop.~3.1]{La22b}).
\begin{thm}\label{thm:kerpro}
 Let $X$, $Y\subseteq {\mathbb{R}}^n$. Let $s_1$, $s_2$, $s_3$, $t_1$, $t_2$, $t_3\in {\mathbb{R}}$. 
\begin{enumerate}
\item[(i)] If $K_1\in  {\mathcal{K}}_{s_{1},s_{2},s_{3}   }(X\times Y)$ and $K_2\in  {\mathcal{K}}_{t_{1},t_{2},t_{3}   }(X\times Y)$, then the following inequality holds
\begin{eqnarray*}
\lefteqn{
|  K_1(x',y)K_2(x',y)- K_1(x'',y)K_2(x'',y)  |
}
\\ \nonumber
&&\qquad\qquad\qquad\qquad \qquad 
\leq\|K_1\|_{ {\mathcal{K}}_{s_{1},s_{2},s_{3}   }(X\times Y)}
\|K_2\|_{ {\mathcal{K}}_{t_{1},t_{2},t_{3}   }(X\times Y)}
\\ \nonumber
&&\qquad\qquad\qquad\qquad \qquad   
\quad
\times\left(
\frac{|x'-x''|^{s_3}}{|x'-y|^{s_2+t_1}}+ \frac{2^{|s_1|}|x'-x''|^{t_3}}{|x'-y|^{t_2+s_1}}
\right)
\end{eqnarray*}
for all $x',x''\in X$, $x'\neq x''$, $y\in Y\setminus{\mathbb{B}}_{n}(x',2|x'-x''|)$. 
\item[(ii)] The pointwise product is bilinear and continuous from
\[
{\mathcal{K}}_{s_{1},s_1+s_3,s_{3}   }(X\times Y)\times  {\mathcal{K}}_{t_{1},t_{1}+s_3,s_{3}   }(X\times Y)
\quad\text{to}\quad
{\mathcal{K}}_{s_1+t_{1},s_{1}+s_3+t_1,s_{3}   }(X\times Y)\,.
\] 
\end{enumerate}
\end{thm}

\begin{prop}\label{prop:prkerho}
 Let $X$, $Y\subseteq {\mathbb{R}}^n$. Let $s_1$, $s_2$, $s_3 \in {\mathbb{R}}$,  
 $\alpha\in]0,1]$. Then the following statements hold.
\begin{enumerate}
\item[(i)] If $K\in {\mathcal{K}}_{s_{1},s_2,s_{3}   }(X\times Y)$ and $f\in C^{0,\alpha}_b(X)$, then
\[
|K(x,y)f(x)|\,|x-y|^{s_1}\leq \|K\|_{ {\mathcal{K}}_{s_{1},X\times Y  } }\sup_X|f|
\quad\forall (x,y)\in X\times Y\setminus{\mathbb{D}}_{X\times Y}\,.
\]
and
\begin{eqnarray*}
\lefteqn{
|K(x',y)f(x')-K(x'',y)f(x'')|
}
\\ \nonumber
&&\quad
\leq \|K\|_{	 {\mathcal{K}}_{s_{1},s_2,s_{3}   }(X\times Y) }\|f\|_{	C^{0,\alpha}_b(X)	}
\left\{
\frac{|x'-x''|^{s_3}}{|x'-y|^{s_2}}+2^{|s_1|}\frac{|x'-x''|^{\alpha}}{|x'-y|^{s_1}} 
\right\}
\end{eqnarray*}
for all $x',x''\in X$, $x'\neq x''$, $y\in Y\setminus{\mathbb{B}}_{n}(x',2|x'-x''|)$. 
\item[(ii)] If  $s_2\geq s_1$ and $X$ and $Y$ are both bounded, then the map from 
\[
{\mathcal{K}}_{s_{1},s_2,s_{3}   }(X\times Y)\times C^{0,s_3}_b(X)\quad\text{to}\quad{\mathcal{K}}_{s_{1},s_2,s_{3}   }(X\times Y)
\]
 that takes the pair $(K,f)$ to the kernel $K(x,y)f(x)$ of the variable $(x,y)\in (X\times Y)\setminus {\mathbb{D}}_{X\times Y}$ is bilinear and continuous. 
 \item[(iii)] The map from 
\[
{\mathcal{K}}_{s_{1},s_2,s_{3}   }(X\times Y)\times C^{0}_b(Y)\quad\text{to}\quad{\mathcal{K}}_{s_{1},s_2,s_{3}   }(X\times Y)
\]
 that takes the pair $(K,f)$ to the kernel $K(x,y)f(y)$ of the variable $(x,y)\in (X\times Y)\setminus {\mathbb{D}}_{X\times Y}$ is bilinear and continuous. 
\end{enumerate}
\end{prop}
Next we have the following imbedding statement that holds for bounded sets (cf.~\cite[Prop.~3.2]{La22b}).
\begin{prop}\label{prop:kerem}
Let $X$, $Y$ be bounded subsets of ${\mathbb{R}}^n$. Let $s_1$, $s_2$, $s_3$, $t_1$, $t_2$, $t_3\in {\mathbb{R}}$. Then the following statements hold.
\begin{enumerate}
\item[(i)] If $t_1\geq s_1$ then $ {\mathcal{K}}_{s_1,X\times Y}$ is continuously embedded into  $ {\mathcal{K}}_{t_1,X\times Y}$.
\item[(ii)] If $t_1\geq s_1$, $t_3\leq s_3$ and $(t_2-t_3)\geq(s_2-s_3)$, then ${\mathcal{K}}_{s_{1},s_{2},s_{3}   }(X\times Y)$ is continuously embedded into ${\mathcal{K}}_{t_{1},t_{2},t_{3}   }(X\times Y)$.
\item[(iii)] If $t_1\geq s_1$,  $t_3\leq s_3$, then ${\mathcal{K}}_{s_{1},s_{1}+s_3,s_{3}   }(X\times Y)$ is continuously embedded into the space ${\mathcal{K}}_{t_{1},t_{1}+t_3,t_{3}   }(X\times Y)$.
\end{enumerate}
\end{prop}
We now show that we can associate a potential type kernel to all H\"{o}lder continuous functions
(cf.~\cite[Lem.~3.3]{La22b}).
\begin{lem}\label{lem:hoker}
 Let $X$, $Y$ be subsets of ${\mathbb{R}}^n$. Let $\alpha\in]0,1]$. Then the following statements hold.
 \begin{enumerate}
\item[(i)] If $\mu\in C^{0,\alpha}(X\cup Y)$, then the map $\Xi[\mu]$  defined by
\begin{equation}\label{lem:hoker1}
\Xi[\mu](x,y)\equiv \mu(x)-\mu(y)\qquad\forall (x,y)\in (X\times Y)\setminus {\mathbb{D}}_{X\times Y}
\end{equation}
 belongs to ${\mathcal{K}}_{-\alpha,0,\alpha  }(X\times Y)$.
\item[(ii)] The operator $\Xi$ from $C^{0,\alpha}(X\cup Y)$ to ${\mathcal{K}}_{-\alpha,0,\alpha  }(X\times Y)$ that takes $\mu$ to $\Xi[\mu]$ is linear and continuous.
\end{enumerate}
 \end{lem}
In order to introduce a result  of \cite[Thm.~6.3]{La22b}, we need to introduce a further norm for kernels
 in the case in which  $Y$ is a compact manifold of class $C^1$ that is imbedded in $M={\mathbb{R}}^n$ and $X=Y$.
\begin{defn}
Let $Y$ be a compact manifold of class $C^1$ that is imbedded in ${\mathbb{R}}^n$.
    Let $s_1$, $s_2$, $s_3\in {\mathbb{R}}$. We set
\begin{eqnarray*}
\lefteqn{
 {\mathcal{K}}_{ s_1, s_2, s_3  }^\sharp(Y\times Y)
  \equiv 
  \biggl\{\biggr.
K\in  {\mathcal{K}}_{ s_1, s_2, s_3  }(Y\times Y):\,
}
\\ \nonumber
&&\ \ 
\sup_{x\in Y}\sup_{r\in ]0,+\infty[}
\left|
\int_{Y\setminus {\mathbb{B}}_n(x,r)}K(x,y)\, d\sigma_y	 
\right|<+\infty
 \biggl.\biggr\}  
\end{eqnarray*}
and
\begin{eqnarray*}
\lefteqn{
\|K\|_{{\mathcal{K}}_{ s_1, s_2, s_3  }^\sharp(Y\times Y)}
\equiv
\|K\|_{{\mathcal{K}}_{ s_1, s_2, s_3  }(Y\times Y)}
}
\\ \nonumber
&&+
\sup_{x\in Y}\sup_{r\in ]0,+\infty[}
\left|
\int_{Y\setminus {\mathbb{B}}_n(x,r)}K(x,y)\,d\sigma_y	 
\right|\quad\forall 
 K\in {\mathcal{K}}_{ s_1, s_2, s_3  }^\sharp(Y\times Y)\,.
 \end{eqnarray*}
\end{defn}
Clearly,  $({\mathcal{K}}^\sharp_{ s_{1},s_{2},s_{3}   }(Y\times Y),\|\cdot\|_{  {\mathcal{K}}^\sharp_{s_{1},s_{2},s_{3}   }(Y\times Y)  })$ is a normed space. By definition, ${\mathcal{K}}^\sharp_{ s_{1},s_{2},s_{3}   }(Y\times Y)$ is continuously embedded into  ${\mathcal{K}}_{ s_{1},s_{2},s_{3}   }(Y\times Y)$. Next we introduce a function that we need for a generalized H\"{o}lder norm. For each $\theta\in]0,1]$, we define the function $\omega_{\theta}(\cdot)$ from $[0,+\infty[$ to itself by setting
\begin{equation}
\label{omth}
\omega_{\theta}(r)\equiv
\left\{
\begin{array}{ll}
0 &r=0\,,
\\
r^{\theta}|\ln r | &r\in]0,r_{\theta}]\,,
\\
r_{\theta}^{\theta}|\ln r_{\theta} | & r\in ]r_{\theta},+\infty[\,,
\end{array}
\right.
\end{equation}
where
$
%\label{omth1}
r_{\theta}\equiv e^{-1/\theta}
$ for all $\theta\in ]0,1]$. Obviously, $\omega_{\theta}(\cdot) $ is concave and satisfies   condition (\ref{om}).
 We also note that if ${\mathbb{D}}\subseteq {\mathbb{R}}^n$, then the continuous embedding
\[
C^{0, \theta }_b({\mathbb{D}})\subseteq 
C^{0,\omega_\theta(\cdot)}_b({\mathbb{D}})\subseteq 
C^{0,\theta'}_b({\mathbb{D}})
\]
holds for all $\theta'\in ]0,\theta[$.  
 We also need to consider convolution kernels, thus we introduce the following notation.
 If $n\in {\mathbb{N}}\setminus\{0\}$, $m\in {\mathbb{N}}$, $h\in {\mathbb{R}}$, $\alpha\in ]0,1]$, then we set 
\begin{equation}\label{volume.klhomog}
{\mathcal{K}}^{m,\alpha}_h \equiv\biggl\{
k\in C^{m,\alpha}_{ {\mathrm{loc}}}({\mathbb{R}}^n\setminus\{0\}):\, k\ {\text{is\ positively\ homogeneous\ of \ degree}}\ h
\biggr\}\,,
\end{equation}
where $C^{m,\alpha}_{ {\mathrm{loc}}}({\mathbb{R}}^n\setminus\{0\})$ denotes the set of functions of 
$C^{m}({\mathbb{R}}^n\setminus\{0\})$ whose restriction to $\overline{\Omega}$ is of class $C^{m,\alpha}(\overline{\Omega})$ for all bounded open subsets $\Omega$ of ${\mathbb{R}}^n$ such that
$\overline{\Omega}\subseteq {\mathbb{R}}^n\setminus\{0\}$
and we set
\[
\|k\|_{ {\mathcal{K}}^{m,\alpha}_h}\equiv \|k\|_{C^{m,\alpha}(\partial{\mathbb{B}}_n(0,1))}\qquad\forall k\in {\mathcal{K}}^{m,\alpha}_h\,.
\]
We can easily verify that $ \left({\mathcal{K}}^{m,\alpha}_h , \|\cdot\|_{ {\mathcal{K}}^{m,\alpha}_h}\right)$ is a Banach space. We also mention the following variant of a well known statement (cf. \textit{e.g.}, \cite[Lem.~3.11]{La22d}).
 \begin{lem}\label{lem:knnnone}
 Let $n\in {\mathbb{N}}\setminus\{0\}$, $h\in [0,+\infty[$. If $k\in C^{0,1}_{ {\mathrm{loc}} }({\mathbb{R}}^n\setminus\{0\})$ is positively homogeneous of degree $-h$, then $k(x-y)\in {\mathcal{K}}_{h,h+1,1}({\mathbb{R}}^n\times {\mathbb{R}}^n)$.
 Moreover, the map from ${\mathcal{K}}^{0,1}_{-h}$   to ${\mathcal{K}}_{h,h+1,1}({\mathbb{R}}^n\times {\mathbb{R}}^n)$ which takes $k$ to $k(x-y)$ is linear and continuous (see (\ref{volume.klhomog}) for the definition of $
{\mathcal{K}}^{0,1}_{-h}$).
\end{lem}

If $X$ and $Y$ are subsets of ${\mathbb{R}}^n$, then the restriction operator  
\[
\text{from}\ {\mathcal{K}}_{h,h+1,1}({\mathbb{R}}^n\times {\mathbb{R}}^n)\ \text{to}\
{\mathcal{K}}_{h,h+1,1}(X\times Y)
\]
is linear and continuous. Thus Lemma \ref{lem:knnnone} implies that if $h\in[0,+\infty[$, then the map
\[
\text{from}\  {\mathcal{K}}^{0,1}_{-h} \ \ \text{ to}\   {\mathcal{K}}_{h,h+1,1}(X\times Y)\,,
 \]
  which takes $k$ to $k(x-y)$ is linear and continuous.

\begin{rem} 
 As Lemma \ref{lem:knnnone} shows the convolution kernels associated to positively homogeneous functions of negative degree are standard kernels. We note however that 
 there exist potential type kernels that belong to a class ${\mathcal{K}}_{s_1,s_2,s_3} (X\times Y)$ with   $s_2\neq s_1+s_3$. 
\end{rem}

\section{Technical preliminaries on the differential operator}\label{sec:techdiop}
  Let $\Omega$ be a bounded open subset of ${\mathbb{R}}^{n}$ of class $C^{1}$.  
 The kernel of the boundary integral operator corresponding to the  double layer potential 
  is the following
\begin{eqnarray}\label{eq:tgdlgen}
\lefteqn{
\overline{B^{*}_{\Omega,y}}\left(S_{{\mathbf{a}}}(x-y)\right)\equiv - \sum_{l,j=1}^{n} a_{jl}\nu_{l}(y)\frac{\partial S_{ {\mathbf{a}} } }{\partial x_{j}}(x-y) 
}
\\ \nonumber
&&\qquad\qquad\qquad
 - \sum_{l=1}^{n}\nu_{l}(y)a_{l}
S_{ {\mathbf{a}} }(x-y)\qquad\forall (x,y)\in (\partial\Omega)^2\setminus {\mathbb{D}}_{\partial\Omega} 
\end{eqnarray}
(cf.~(\ref{introd3})). In order to analyze the kernel of the double layer potential, we need some more information on the fundamental solution $S_{ {\mathbf{a}} } $. To do so, we introduce the fundamental solution $S_{n}$ of the Laplace operator. Namely, we set
\[
%\label{ups}
S_{n}(x)\equiv
\left\{
\begin{array}{lll}
\frac{1}{s_{n}}\ln  |x| \qquad &   \forall x\in 
{\mathbb{R}}^{n}\setminus\{0\},\quad & {\mathrm{if}}\ n=2\,,
\\
\frac{1}{(2-n)s_{n}}|x|^{2-n}\qquad &   \forall x\in 
{\mathbb{R}}^{n}\setminus\{0\},\quad & {\mathrm{if}}\ n>2\,,
\end{array}
\right.
\]
where $s_{n}$ denotes the $(n-1)$ dimensional measure of 
$\partial{\mathbb{B}}_{n}(0,1)$ and
we follow a formulation of Dalla Riva \cite[Thm.~5.2, 5.3]{Da13} and Dalla Riva, Morais and Musolino \cite[Thm.~5.5]{DaMoMu13}, that we state  as in Dondi and the author  \cite[Cor.~4.2]{DoLa17}  (see also John~\cite{Jo55}, and Miranda~\cite{Mi65} for homogeneous operators, and Mitrea and Mitrea~\cite[p.~203]{MitMit13}).   
 \begin{prop}
 \label{prop:ourfs} 
Let ${\mathbf{a}}$ be as in (\ref{introd0}), (\ref{ellip}), (\ref{symr}). 
Let $S_{ {\mathbf{a}} }$ be a fundamental solution of $P[{\mathbf{a}},D]$. 
Then there exist an invertible matrix $T\in M_{n}({\mathbb{R}})$ such that
\begin{equation}
\label{prop:ourfs0}
a^{(2)}=TT^{t}\,,
\end{equation}
 a real analytic function $A_{1}$ from $\partial{\mathbb{B}}_{n}(0,1)\times{\mathbb{R}}$ to ${\mathbb{C}}$ such that
 $A_{1}(\cdot,0)$ is odd,    $b_{0}\in {\mathbb{C}}$, a real analytic function $B_{1}$ from ${\mathbb{R}}^{n}$ to ${\mathbb{C}}$ such that $B_{1}(0)=0$, and a real analytic function $C $ from ${\mathbb{R}}^{n}$ to ${\mathbb{C}}$ such that
\begin{eqnarray}%attenzione a formulare diversamente nel lavoro
\label{prop:ourfs1}
\lefteqn{
S_{ {\mathbf{a}} }(x)
= 
\frac{1}{\sqrt{\det a^{(2)} }}S_{n}(T^{-1}x)
}
\\ \nonumber
&&\qquad
+|x|^{3-n}A_{1}(\frac{x}{|x|},|x|)
 +(B_{1}(x)+b_{0}(1-\delta_{2,n}))\ln  |x|+C(x)\,,
 %\\ \nonumber
%&&\qquad
%-\frac{\delta_{2,n}}{2\pi}\int_{\partial{\mathbb{B}}_{2}(0,1)}S_{a^{(2)}}-A_{0}\,d\sigma
\end{eqnarray}
for all $x\in {\mathbb{R}}^{n}\setminus\{0\}$,
 and such that both $b_{0}$ and $B_{1}$   equal zero
if $n$ is odd. Moreover, 
 \[
 \frac{1}{\sqrt{\det a^{(2)} }}S_{n}(T^{-1}x) 
 \]
is a fundamental solution for the principal part
  of $P[{\mathbf{a}},D]$.
\end{prop}
In particular for the statement that $A_{1}(\cdot,0)$ is odd, we refer to
Dalla Riva, Morais and Musolino \cite[Thm.~5.5, (32)]{DaMoMu13}, where $A_{1}(\cdot,0)$ coincides with ${\mathbf{f}}_1({\mathbf{a}},\cdot)$ in that paper. Here we note that a function $A$ from $(\partial{\mathbb{B}}_{n}(0,1))\times{\mathbb{R}}$ to ${\mathbb{C}}$ is said to be real analytic provided that it has a real analytic extension   to an open neighbourhood of $(\partial{\mathbb{B}}_{n}(0,1))\times{\mathbb{R}}$ in 
${\mathbb{R}}^{n+1}$. Then we have the following elementary lemma (cf.~\textit{e.g.}, \cite[Lem.~4.2]{La22d}).
\begin{lem}\label{lem:anexsph}
 Let $n\in {\mathbb{N}}\setminus\{0,1\}$. A function $A$ from   $(\partial{\mathbb{B}}_{n}(0,1))\times{\mathbb{R}}$ to ${\mathbb{C}}$ is  real analytic if and only if the function $\tilde{A}$ from $({\mathbb{R}}^n\setminus\{0\}) \times{\mathbb{R}}$ defined by 
\begin{equation}\label{lem:anexsph1}
\tilde{A}(x,r)\equiv A(\frac{x}{|x|},r)\qquad\forall (x,r)\in ({\mathbb{R}}^n\setminus\{0\}) \times{\mathbb{R}}
\end{equation}
is real analytic.
 \end{lem}
Then  one can prove the following formula for the gradient of the fundamental solution (see Dondi and the author  \cite[Lem.~4.3,  (4.8) and the following 2 lines]{DoLa17}. Here one should remember that $A_1(\cdot,0)$ is odd and that $b_0=0$ if $n$ is odd). 
\begin{prop}
\label{prop:grafun}
 Let ${\mathbf{a}}$ be as in (\ref{introd0}), (\ref{ellip}), (\ref{symr}). Let $T\in M_{n}({\mathbb{R}})$  be as in (\ref{prop:ourfs0}). Let $S_{ {\mathbf{a}} }$ be a fundamental solution of $P[{\mathbf{a}},D]$. Let  $B_{1}$, $C$
 be as in Proposition \ref{prop:ourfs}. 
  Then there exists a real analytic function $A_{2}$ from $\partial{\mathbb{B}}_{n}(0,1)\times{\mathbb{R}}$ to ${\mathbb{C}}^{n}$ such that
\begin{eqnarray}
\label{prop:grafun1}
\lefteqn{
DS_{ {\mathbf{a}} }(x)=\frac{1}{ s_{n}\sqrt{\det a^{(2)} } }
|T^{-1}x|^{-n}x^{t}(a^{(2)})^{-1}
}
\\ \nonumber
&& 
+|x|^{2-n}A_{2}(\frac{x}{|x|},|x|)+DB_{1}(x)\ln |x|+DC(x)
\quad\forall x\in {\mathbb{R}}^{n}\setminus\{0\}\,.
\end{eqnarray}
Moreover,   $A_2(\cdot,0)$ is even.
\end{prop}
Next we introduce the following technical lemma 
(see Dondi and the author  \cite[Lem. 3.2 (v), 3.3]{DoLa17}). See also \cite[Lem.~4.5]{La22d}.
\begin{lem}
\label{lem:fanes}
Let $Y$ be a nonempty bounded subset of ${\mathbb{R}}^{n}$. Then the following statements hold.  
\begin{enumerate}
\item[(i)] Let ${\mathrm{diam}}\,(Y)$ be the diameter of $Y$, $F\in {\mathrm{Lip}}(\partial{\mathbb{B}}_{n}(0,1)\times [0,{\mathrm{diam}}\,(Y)])$ with
\begin{eqnarray*}
\lefteqn{
{\mathrm{Lip}}(F)
\equiv\biggl\{\biggr.
\frac{|F(\theta',r')-F(\theta'',r'')|}{ |\theta'-\theta''|+|r'-r''| }:\,
}
\\ \nonumber
&&\quad 
(\theta',r'),(\theta'',r'')\in \partial{\mathbb{B}}_{n}(0,1)\times [0,{\mathrm{diam}}\,(Y)],\ (\theta',r')\neq (\theta'',r'')
\biggl.\biggr\}\,.
\end{eqnarray*}
Then 
\begin{eqnarray}
\label{lem:fanes1}
\lefteqn{
\left|
F\left(
\frac{x'-y}{|x'-y|},|x'-y|
\right)
-
F\left(
\frac{x''-y}{|x''-y|},|x''-y|
\right)
\right|
}
\\ \nonumber
&&\quad
\leq
{\mathrm{Lip}}(F) (2+   {\mathrm{diam}}\,(Y))
\frac{|x'-x''|}{|x'-y|}\,
\quad\forall y\in Y
\setminus {\mathbb{B}}_{n}(x',2|x'-x''|)\,,
\end{eqnarray}
for all $x',x''\in Y$, $x'\neq x''$. In particular, if $f\in C^{1}(\partial{\mathbb{B}}_{n}(0,1)\times{\mathbb{R}},{\mathbb{C}})$, then
\begin{eqnarray*}
\lefteqn{
M_{f,Y}\equiv\sup\biggl\{\biggr.
\left|
f\left(
\frac{x'-y}{|x'-y|},|x'-y|
\right)
-
f\left(
\frac{x''-y}{|x''-y|},|x''-y|
\right)
\right| 
}
\\ \nonumber
&& \qquad 
\times\frac{|x'-y|}{|x'-x''|}:\,x',x''\in Y, x'\neq x'',  y\in Y   
\setminus {\mathbb{B}}_{n}(x',2|x'-x''|)
\biggl.\biggr\} 
\end{eqnarray*}
is finite and thus the kernel $f\left(
\frac{x-y}{|x-y|},|x-y|
\right)$ belongs to  ${\mathcal{K}}_{0,1,1}(Y\times Y)$.

\item[(ii)] Let $W$ be an open neighbourhood of $\overline{Y-Y}$. Let $f\in C^{1}(W,{\mathbb{C}})$. Then 
\begin{eqnarray*}
\lefteqn{
\tilde{M}_{f,Y}\equiv
 \sup\biggl\{\biggr.
|
f(x'-y)-f(x''-y)|\,|x'-x''|^{-1}
}
\\ \nonumber
&& \qquad\qquad\qquad\qquad\qquad\quad
:\,x',x''\in Y, x'\neq x'',  y\in Y 
\biggl.\biggr\}<+\infty\,.
\end{eqnarray*}
Here $Y-Y\equiv\{y_{1}-y_{2}:\ y_{1}, y_{2}\in Y\}$. In particular, the kernel $f(x-y)$ belongs to the class ${\mathcal{K}}_{0,0,1}(Y\times Y)$, which is continuously imbedded into ${\mathcal{K}}_{0,1,1}(Y\times Y)$.
\item[(iii)] The kernel $\ln|x-y|$ belongs to ${\mathcal{K}}_{\epsilon,1,1}(Y\times Y)$ for all $\epsilon\in]0,1[$.
\end{enumerate}
\end{lem}

In order to prove regularity results for the double layer potential, we need the definition of tangential derivative and some auxiliary operator that we now introduce.
Let $\Omega$ be an open subset of ${\mathbb{R}}^n$ of class $C^1$.  If $l,r\in\{1,\dots,n\}$,  then $M_{lr}$ denotes the tangential derivative 
 operator from $C^{1}(\partial\Omega)$ to $C^{0}(\partial\Omega)$ that takes $f$ to  
 \begin{equation}
\label{mlr}
M_{lr}[f]\equiv \nu_{l}\frac{\partial\tilde{f}}{\partial x_{r}}-
\nu_{r}\frac{\partial\tilde{f}}{\partial x_{l}}\qquad {\text{on}}\ \partial\Omega\,,
\end{equation}
where  $\tilde{f}$ is any continuously differentiable  extension of $f$  to an open neighborhood of $\partial\Omega$. We note that $M_{lr}[f]$ is independent of the specific choice of $\tilde{f}$ (cf.~\textit{e.g.}, Dalla Riva, the author and Musolino \cite[\S 2.21]{DaLaMu21}). For the definition of tangential gradient $ {\mathrm{grad}}_{\partial\Omega}$, we refer   to Kirsch and Hettlich \cite[A.5]{KiHe15}, Chavel~\cite[Chap.~1]{Cha84}.  Then we set
\begin{equation}\label{qrs0}
Q_j[g,\mu](x)
=\int_{\partial\Omega}(g(x)-g(y))\frac{\partial S_{ {\mathbf{a}} }}{\partial x_{j}}(x-y)\mu(y)\,d\sigma_{y}\quad\forall x\in \partial\Omega\,,
\end{equation}
for all $(g,\mu)\in  C^{0,1}(\partial\Omega)\times L^{\infty}(\partial\Omega)$ for all $j\in\{1,\dots,n\}$. As a first step, we prove the following technical statement that determines the second order partial derivatives of the kernel $S_{ {\mathbf{a}} }(x-y)$, the class membership of the corresponding kernels and the class
of the tangential gradient of the kernel $\frac{\partial S_{ {\mathbf{a}} }}{\partial x_{j}}(x-y)$ with respect to its first variable on the boundary of an open set of class $C^{1,\alpha}$
 for all $j\in\{1,\dots,n\}$.
\begin{lem}\label{lem:tgdjsgen}  
 Let ${\mathbf{a}}$ be as in (\ref{introd0}), (\ref{ellip}), (\ref{symr}).  Let $S_{ {\mathbf{a}} }$ be a fundamental solution of $P[{\mathbf{a}},D]$.    Let $j,h\in\{1,\dots,n\}$.   Then the following statements hold.
 \begin{enumerate}
\item[(i)]   
\begin{eqnarray}\label{lem:tgdjsgen1}  
\lefteqn{
\frac{\partial}{\partial x_h}\frac{\partial S_{{\mathbf{a}}}}{\partial x_j}(x-y) 
}
\\ \nonumber
&& 
=\frac{-n
|T^{-1}(x-y)|^{-n-1}
}{s_{n}\sqrt{\det a^{(2)} }} \frac{
\sum_{s,t=1}^n(T^{-1})_{st}(x_t-y_t)(T^{-1})_{sh}
}{|T^{-1}(x-y)|}
\\ \nonumber
&& \quad
\times \sum_{s=1}^n(x_s-y_s)((a^{(2)})^{-1})_{sj}
+
\frac{|T^{-1}(x-y)|^{-n}
}{s_{n}\sqrt{\det a^{(2)} }} ((a^{(2)})^{-1})_{hj}
\\ \nonumber
&& \quad
+(2-n)|x-y|^{1-n}\frac{x_h-y_h}{|x-y|}A_{2,j}\left(\frac{x-y}{|x-y|},|x-y|\right)
\\ \nonumber
&& \quad
+|x-y|^{2-n}\biggl\{\sum_{s=1}^n\frac{\partial A_{2,j}}{\partial x_s}\left(\frac{x-y}{|x-y|},|x-y|\right)
\\ \nonumber
&& \quad
\times \biggl(
\delta_{sh}|x-y|-\frac{(x_s-y_s)(x_h-y_h)}{|x-y|}
\biggr)|x-y|^{-2}
\\ \nonumber
&& \quad
+\frac{\partial A_{2,j}}{\partial r}
\left(\frac{x-y}{|x-y|},|x-y|\right)\frac{x_h-y_h}{|x-y|}\biggr\}
\\ \nonumber
&& \quad
+\frac{\partial^2B_1}{\partial x_h\partial x_j}(x-y)\ln|x-y|
+\frac{\partial B_1}{\partial x_j}(x-y)\frac{x_h-y_h}{|x-y|^2}
\\ \nonumber
&& \quad
+\frac{\partial^2C}{\partial x_h\partial x_j}(x-y)
\end{eqnarray}
for all $x$, $y\in {\mathbb{R}}^n$, $x\neq y$. 
\item[(ii)] If $G$ be a nonempty bounded subset of ${\mathbb{R}}^{n}$, then the kernel $\frac{\partial}{\partial x_h}\frac{\partial S_{{\mathbf{a}}}}{\partial x_j}(x-y)$ belongs to ${\mathcal{K}}_{n,n+1,1}(G\times G)$. 
\item[(iii)] Let $\alpha\in]0,1]$.  Let $\Omega$ be a bounded open subset of ${\mathbb{R}}^{n}$ of class $C^{1,\alpha}$. Then
\begin{eqnarray}\label{lem:tgdjsgen2}
\lefteqn{
\left({\mathrm{grad}}_{\partial\Omega,x} \left(\frac{\partial S_{{\mathbf{a}}}}{\partial x_j}(x-y)\right) \right)_h
}
\\ \nonumber
&&\qquad
=\frac{\partial}{\partial x_h}\frac{\partial S_{{\mathbf{a}}}}{\partial x_j}(x-y)
-\sum_{l=1}^{n}\frac{\partial}{\partial x_l}\frac{\partial S_{{\mathbf{a}}}}{\partial x_j}(x-y)\nu_l(x)\nu_h(x)\,,
 \end{eqnarray}
 for all $(x,y)\in (\partial\Omega)^2\setminus {\mathbb{D}}_{\partial\Omega}$ and the kernel $\left({\mathrm{grad}}_{\partial\Omega,x} \left(\frac{\partial S_{{\mathbf{a}}}}{\partial x_j}(x-y) \right)\right)_h$ belongs to ${\mathcal{K}}_{n,n+ \alpha,\alpha}((\partial\Omega)\times (\partial\Omega))$. Here ${\mathrm{grad}}_{\partial\Omega,x}$ denotes the tangential gradient with respect to the $x$ variable. 
 \end{enumerate}
\end{lem}
{\bf Proof.} Statement (i) holds by formula (\ref{prop:grafun1}) and by standard differentiation rules. (ii) By Lemma \ref{lem:knnnone}, the kernel  
\[
\frac{-n
|T^{-1}(x-y)|^{-n-1}
}{s_{n}\sqrt{\det a^{(2)} }} \frac{
\sum_{s,t=1}^n(T^{-1})_{st}(x_t-y_t)(T^{-1})_{sh}
}{|T^{-1}(x-y)|} \sum_{s=1}^n(x_s-y_s)((a^{(2)})^{-1})_{sj}
\]
and the kernel
\[
\frac{|T^{-1}(x-y)|^{-n}
}{s_{n}\sqrt{\det a^{(2)} }} ((a^{(2)})^{-1})_{hj}
\]
belong to ${\mathcal{K}}_{n,n+1,1}(G\times G)$. Since $A_{2,j}$ is real analytic in $\partial{\mathbb{B}}_n(0,1)\times {\mathbb{R}}$, Lemma \ref{lem:fanes} (i) implies that 
$A_{2,j}\left(\frac{x-y}{|x-y|},|x-y|\right)$ belongs to ${\mathcal{K}}_{0,1,1}(G\times G)$. By Lemma \ref{lem:knnnone}, the kernel $|x-y|^{1-n}\frac{x_h-y_h}{|x-y|}$ belongs to
${\mathcal{K}}_{n-1,n,1}(G\times G)$. Then the product Theorem \ref{thm:kerpro} (ii) implies that the product
\[
(2-n)|x-y|^{1-n}\frac{x_h-y_h}{|x-y|}A_{2,j}(\frac{x-y}{|x-y|},|x-y|)
\]
belongs to ${\mathcal{K}}_{n-1,n,1}(G\times G)$.

Since $\frac{\partial A_{2,j}}{\partial x_s}$ is real analytic in $\partial{\mathbb{B}}_n(0,1)\times {\mathbb{R}}$, Lemma \ref{lem:fanes} (i) implies that 
$\frac{\partial A_{2,j}}{\partial x_s}\left(\frac{x-y}{|x-y|},|x-y|\right)$  belongs to ${\mathcal{K}}_{0,1,1}(G\times G)$. 
By Lemma \ref{lem:knnnone}, the kernel 
\[
|x-y|^{2-n}\biggl(
\delta_{sh}|x-y|-\frac{(x_s-y_s)(x_h-y_h)}{|x-y|}
\biggr)|x-y|^{-2}
\]
belongs to ${\mathcal{K}}_{n-1,n,1}(G\times G)$. Then the product Theorem \ref{thm:kerpro} (ii) implies that the product
\[
|x-y|^{-n} \sum_{s=1}^n\frac{\partial A_{2,j}}{\partial x_s}\left(\frac{x-y}{|x-y|},|x-y|\right)
  \biggl(
\delta_{sh}|x-y|-\frac{(x_s-y_s)(x_h-y_h)}{|x-y|}
\biggr) 
\]
belongs to ${\mathcal{K}}_{n-1,n,1}(G\times G)$.

Since $\frac{\partial A_{2,j}}{\partial r}$ is real analytic in $\partial{\mathbb{B}}_n(0,1)\times {\mathbb{R}}$, Lemma \ref{lem:fanes} (i) implies that 
$\frac{\partial A_{2,j}}{\partial r}\left(\frac{x-y}{|x-y|},|x-y|\right)$  belongs to ${\mathcal{K}}_{0,1,1}(G\times G)$. 
By Lemma \ref{lem:knnnone}, the kernel 
\[
|x-y|^{2-n} \frac{x_h-y_h}{|x-y|}
\]
belongs to ${\mathcal{K}}_{n-2,n-1,1}(G\times G)$. Then the product Theorem \ref{thm:kerpro} (ii) implies that the product
\[
|x-y|^{2-n}\frac{\partial A_2}{\partial r}
\left(\frac{x-y}{|x-y|},|x-y|\right)\frac{x_h-y_h}{|x-y|}
\]
belongs to ${\mathcal{K}}_{n-2,n-1,1}(G\times G)$. 

Since $B_1$ is real analytic, Lemma \ref{lem:fanes} (ii) implies that  
$
\frac{\partial^2B_1}{\partial x_h\partial x_j}(x-y) 
$
belongs to ${\mathcal{K}}_{0,1,1}(G\times G)$. By Lemma \ref{lem:fanes} (iii), the kernel  $\ln|x-y|$ belongs to ${\mathcal{K}}_{\epsilon,1,1}(G\times G)$ for all $\epsilon\in]0,1[$. By the embedding Proposition \ref{prop:kerem} (ii), ${\mathcal{K}}_{\epsilon,1,1}(G\times G)$ is contained in 
${\mathcal{K}}_{\epsilon,\epsilon+1,1}(G\times G)$  for all $\epsilon\in]0,1[$.

Then the product Theorem \ref{thm:kerpro} (ii) implies that the product
\[
\frac{\partial^2B_1}{\partial x_h\partial x_j}(x-y)\ln|x-y|
\]
belongs to ${\mathcal{K}}_{\epsilon,\epsilon+1,1}(G\times G)$  for all $\epsilon\in]0,1[$.

Since $B_1$ is real analytic, Lemma \ref{lem:fanes} (ii) implies that the kernel
$
\frac{\partial B_1}{\partial x_j}(x-y)
$
belongs to ${\mathcal{K}}_{0,1,1}(G\times G)$. By Lemma \ref{lem:knnnone}, the kernel 
$
\frac{x_h-y_h}{|x-y|^2}
$
belongs to ${\mathcal{K}}_{1,2,1}(G\times G)$. Then the product Theorem \ref{thm:kerpro} (ii) implies that the product
\[
\frac{\partial B_1}{\partial x_j}(x-y)\frac{x_h-y_h}{|x-y|^2}
\]
belongs to ${\mathcal{K}}_{1,2,1}(G\times G)$.

Since $C$ is real analytic, Lemma \ref{lem:fanes} (ii) implies that  
$
\frac{\partial^2C}{\partial x_h\partial x_j}(x-y)
$
belongs to ${\mathcal{K}}_{0,1,1}(G\times G)$. Thus we have proved that each addendum   in the right hand side of formula (\ref{lem:tgdjsgen1})  is contained in one of the following classes  
\begin{eqnarray*}
&&{\mathcal{K}}_{n,n+1,1}(G\times G)\,,\quad
{\mathcal{K}}_{n-1,n,1}(G\times G)\,,\quad
{\mathcal{K}}_{n-2,n-1,1}(G\times G)\,,\quad
\\ \nonumber
&&
{\mathcal{K}}_{\epsilon,\epsilon+1,1}(G\times G)\quad\forall \epsilon\in]0,1[\,,\quad
 {\mathcal{K}}_{1,2,1}(G\times G)\,,\qquad
{\mathcal{K}}_{0,1,1}(G\times G) \,.
\end{eqnarray*}
Now the imbedding Proposition \ref{prop:kerem} (iii) implies that each of such classes 
is contained in ${\mathcal{K}}_{n,n+1,1}(G\times G)$ and thus the proof of statement (ii) is complete.\par 

(iii) Formula (\ref{lem:tgdjsgen2}) holds by the definition of tangential gradient. By the elementary Lemma \ref{lem:kelem}, we have
\begin{eqnarray*}
\lefteqn{
{\mathcal{K}}_{n,n+1,1}((\partial\Omega)\times (\partial\Omega))
}
\\ \nonumber
&&\qquad
\subseteq {\mathcal{K}}_{n,n+1-(1-\alpha),1-(1-\alpha)}((\partial\Omega)\times (\partial\Omega))
={\mathcal{K}}_{n,n+\alpha,\alpha}((\partial\Omega)\times (\partial\Omega))\,.
\end{eqnarray*}
Then the membership of the components of $\nu $ in $C^{0,\alpha}(\partial\Omega)$,  statement (i) with $G=\partial\Omega$ condition $n\leq n+\alpha$ and the product Proposition \ref{prop:prkerho} (ii) imply that the right hand side of formula (\ref{lem:tgdjsgen2}) defines a kernel of class ${\mathcal{K}}_{n,n+\alpha,\alpha}((\partial\Omega)\times (\partial\Omega))$ and thus the proof is complete.\hfill  $\Box$ 

\vspace{\baselineskip}

Then we introduce the following technical statement 
(cf.~\cite[Thm.~3.2]{La23b}).  
\begin{thm}\label{lem:maolg}
 Let $n\in {\mathbb{N}}$,  $n\geq 2$,   $\tilde{\alpha}\in]0,1]$.  Let $\Omega$ be a bounded open subset of ${\mathbb{R}}^n$ of class $C^{1,\tilde{\alpha}}$. Then there exists   $c^*_{\partial\Omega,\tilde{\alpha}}\in]0,+\infty[$ such that 
\begin{equation}\label{lem:maolg1}
\sup_{x\in \partial\Omega}
\sup_{\epsilon\in ]0,+\infty[}
\left|
\int_{(\partial\Omega)\setminus{\mathbb{B}}_{n}(x,\epsilon)}
k(x-y)\,d\sigma_y
\right|
\leq c^*_{\partial\Omega,\tilde{\alpha}}\left\| k \right\|_{{\mathcal{K}}^{0,1}_{-(n-1) }}
\ \ \forall k\in {\mathcal{K}}^{0,1}_{-(n-1);o }\,,
\end{equation}
where
$
{\mathcal{K}}^{0,1}_{-(n-1);o } \equiv \{k\in {\mathcal{K}}^{0,1}_{-(n-1)}:\, k\ \text{is\ odd}\}
$ (cf.~(\ref{volume.klhomog})).
\end{thm}
 Next we prove the following technical lemma.
\begin{lem}\label{lem:mahgnfs}
 Let $n\in {\mathbb{N}}\setminus\{0,1\}$. 
 Let ${\mathbf{a}}$ be as in (\ref{introd0}), (\ref{ellip}), (\ref{symr}).  Let $S_{ {\mathbf{a}} }$ be a fundamental solution of $P[{\mathbf{a}},D]$.   Let $\alpha\in]0,1]$.  Let $\Omega$ be a bounded open subset of ${\mathbb{R}}^n$  of class $C^{1, \alpha }$.
 Let $j,h,z\in \{1,\dots,n\}$.Then
 \begin{equation}\label{lem:mahgnfs1}
\sup_{x\in \partial\Omega}
\sup_{r\in ]0,+\infty[}
\left|\,
\int_{(\partial\Omega)\setminus{\mathbb{B}}_{n}(x,r)}
(x_z-y_z)\frac{\partial}{\partial x_h}\frac{\partial S_{{\mathbf{a}}}}{\partial x_j}(x-y) 
\,d\sigma_y\right|
<+\infty\,.
\end{equation}
\end{lem}
{\bf Proof.} By formula (\ref{lem:tgdjsgen1}) and by the known inequalities 
\begin{equation}\label{thm:maggnfs2}
\sup_{x\in \partial\Omega}\int_{\partial\Omega}|x-y|^{-\gamma}\,d\sigma_y<+\infty\,,
\quad
\sup_{x\in \partial\Omega}\int_{\partial\Omega}
|\ln |x-y||
\,d\sigma_y<+\infty
\end{equation}
for $\gamma\in]-\infty,(n-1)[$ (cf.~\textit{e.g.},  \cite[Lem.~3.5]{DoLa17}),  we have
\begin{eqnarray}\label{lem:mahgnfs2}
\lefteqn{
\sup_{x\in \partial\Omega}
\sup_{r\in ]0,+\infty[}
\biggl| \,
\int_{(\partial\Omega)\setminus{\mathbb{B}}_{n}(x,r)}
(x_z-y_z)\frac{\partial}{\partial x_h}\frac{\partial S_{{\mathbf{a}}}}{\partial x_j}(x-y) 
}
\\ \nonumber
&&\ \ 
-(x_z-y_z)
\frac{-n
|T^{-1}(x-y)|^{-n-1}
}{s_{n}\sqrt{\det a^{(2)} }} \frac{
\sum_{s,t=1}^n(T^{-1})_{st}(x_t-y_t)(T^{-1})_{sh}
}{|T^{-1}(x-y)|}
\\ \nonumber
&& \ \ 
\times \sum_{s=1}^n(x_s-y_s)((a^{(2)})^{-1})_{sj}
-(x_z-y_z)
\frac{|T^{-1}(x-y)|^{-n}
}{s_{n}\sqrt{\det a^{(2)} }} ((a^{(2)})^{-1})_{hj}
\,d\sigma_y
\biggr|
\\ \nonumber
&& 
\leq 
\sup_{x\in \partial\Omega}
\int_{\partial\Omega}|2-n||x-y|^{2-n}
\left|A_{2,j}\left(\frac{x-y}{|x-y|},|x-y|\right)\right|
\\ \nonumber
&&  \ \ 
+|x-y|^{3-n}\biggl\{\sum_{s=1}^n\left|\frac{\partial A_{2,j}}{\partial x_s}\left(\frac{x-y}{|x-y|},|x-y|\right)\right|2|x-y|^{-1}
\\ \nonumber
&& \ \ 
+\left|\frac{\partial A_{2,j}}{\partial r}
\left(\frac{x-y}{|x-y|},|x-y|\right)\right| \biggr\}
+\left|\frac{\partial^2B_1}{\partial x_h\partial x_j}(x-y)\right||x-y|\ln|x-y|
\\ \nonumber
&& \ \ 
+\left|\frac{\partial B_1}{\partial x_j}(x-y)\right|
+\left|\frac{\partial^2C}{\partial x_h\partial x_j}(x-y)\right||x-y|\,d\sigma_y<+\infty\,.
\end{eqnarray}
Since the function 
\begin{eqnarray*}
\lefteqn{
 \xi_z
\frac{-n
|T^{-1}\xi|^{-n-1}
}{s_{n}\sqrt{\det a^{(2)} }} \frac{
\sum_{s,t=1}^n(T^{-1})_{st}\xi_t(T^{-1})_{sh}
}{|T^{-1}\xi|}
\sum_{s=1}^n\xi_s((a^{(2)})^{-1})_{sj}
}
\\ \nonumber
&& \qquad\qquad
  + \xi_z
\frac{|T^{-1}\xi|^{-n}
}{s_{n}\sqrt{\det a^{(2)} }} ((a^{(2)})^{-1})_{hj}
\qquad\forall\xi\in{\mathbb{R}}^n\setminus\{0\}\end{eqnarray*} 
 is positively homogeneous of degree $-(n-1)$, Theorem \ref{lem:maolg} implies  that 
\begin{eqnarray*}
\lefteqn{\sup_{x\in \partial\Omega}
\sup_{r\in ]0,+\infty[}
\biggl|\,
\int_{(\partial\Omega)\setminus{\mathbb{B}}_{n}(x,r)}
 (x_z-y_z)
\frac{-n
|T^{-1}(x-y)|^{-n-1}
}{s_{n}\sqrt{\det a^{(2)} }} 
}
\\ \nonumber
&& \ \ 
\times 
\frac{
\sum_{s,t=1}^n(T^{-1})_{st}(x_t-y_t)(T^{-1})_{sh}
}{|T^{-1}(x-y)|}
\sum_{s=1}^n(x_s-y_s)((a^{(2)})^{-1})_{sj}
\\ \nonumber
&& \ \ 
+(x_z-y_z)
\frac{|T^{-1}(x-y)|^{-n}
}{s_{n}\sqrt{\det a^{(2)} }} ((a^{(2)})^{-1})_{hj}
\biggr|
\end{eqnarray*}
is finite. Then the above  inequality implies the validity of the statement.
\hfill  $\Box$ 

\vspace{\baselineskip}

\section{An extension of a classical theorem for the single layer potential}
We plan to prove the following extension of a known 
classical result for the single layer potential 
\begin{equation}\label{eq:silapo}
v_\Omega[S_{ {\mathbf{a}} },\mu](x)\equiv 
\int_{\partial\Omega}S_{ {\mathbf{a}} }(x-y)\mu(y)\,d\sigma_{y}
\qquad\forall x\in{\mathbb{R}}^{n}\,,
\end{equation}
for all $\mu\in C^{0,\alpha}(\partial\Omega)$ 
(cf. Miranda~\cite{Mi65}, Kirsch \cite[Thm.~3.3 (a)]{Ki89}, 
Wiegner~\cite{Wi93}, Dalla Riva \cite{Da13}, Dalla Riva, Morais and Musolino \cite{DaMoMu13} and references therein.) 
\begin{thm}
\label{thm:slay}
Let ${\mathbf{a}}$ be as in (\ref{introd0}), (\ref{ellip}), (\ref{symr}). Let $S_{ {\mathbf{a}} }$ be a fundamental solution of $P[{\mathbf{a}},D]$. 
Let   $m\in {\mathbb{N}}\setminus\{0\}$. Let $\Omega$ be a bounded open subset of ${\mathbb{R}}^{n}$ of class $C^{m,1}$. If $\mu\in C^{m-1,1}(\partial\Omega)$, then the restriction
\begin{equation}\label{eq:silapobd}
V_\Omega[S_{{\mathbf{a}}} ,\mu] \equiv v_\Omega[S_{ {\mathbf{a}} },\mu]_{|\partial\Omega}
\end{equation}
 belongs to $C^{m,\omega_1(\cdot)}(\partial\Omega)$. Moreover, the map from the space  $C^{m-1,1}(\partial\Omega)$ to $C^{m,\omega_1(\cdot)}(\partial\Omega)$   that takes $\mu$  to $V_\Omega[S_{{\mathbf{a}}} ,\mu] $ is continuous.
\end{thm}
{\bf Proof.} We proceed by induction on $m$. Let  $m=1$. By the definition of norm in $C^{1,\omega_1(\cdot)}(\partial\Omega)$, it suffices to show that
\begin{enumerate}
\item[(j)] $V_\Omega[S_{{\mathbf{a}}} ,\cdot]$ is linear and continuous from  $C^{0,1}(\partial\Omega)$ to $C^{0}(\partial\Omega)$.
\item[(jj)] $V_\Omega[S_{{\mathbf{a}}} ,\mu]$ is continuously differentiable on $\partial\Omega$ for all $\mu$ in $C^{0,1}(\partial\Omega)$.
\item[(jjj)] $M_{jl}[V_\Omega[S_{{\mathbf{a}}} ,\cdot]]$ is linear and continuous from  $C^{0,1}(\partial\Omega)$ to $C^{0,\omega_1(\cdot)}(\partial\Omega)$ for all $j,l\in \{1,\dots,n\}$,
 \end{enumerate}
 (cf.~\textit{e.g.}, \cite[Lem.~2.3]{DoLa17}).  Since $C^{1,1}(\partial\Omega)$ is continuously imbedded into $C^{1,\alpha}(\partial\Omega)$ for all $\alpha\in]0,1[$, statements (j), (jj) hold by \cite[Thm.~7.1 (i)]{DoLa17}. 
 
 We now consider statement (jjj). To do so, we plan to write a formula for the tangential derivatives of the single layer in terms of the tangential derivatives of the density on $\partial\Omega$.  Let $v^+_\Omega[S_{ {\mathbf{a}} },\mu] $ denote the restriction of
 $v_\Omega[S_{ {\mathbf{a}} },\mu] $ to $\overline{\Omega}$. 
 Since $\Omega$ is of class $C^{1,1}$ and accordingly of class $C^{1,\alpha}$ for all $\alpha\in]0,1[$ and $\mu\in C^{0,\alpha}(\partial\Omega)$, we know that $v^+_\Omega[S_{ {\mathbf{a}} },\mu]\in C^{1,\alpha}(\overline{\Omega})$  for all $\alpha\in]0,1[$ (cf.~\cite[Thm.~7.1 (i)]{DoLa17}).  Let $j,l\in \{1,\dots,n\}$. 
 Since $\Omega$ is of class $C^{1,1}$ and $\nu$ is of class $C^{0,1}$, there exists $\tilde{\nu}\in C^{0,1}({\mathbb{R}}^n)$ with compact support such that
 $\tilde{\nu}_{|\partial\Omega}=\nu$ (cf.~\textit{e.g.}, \cite[Thm.~2.85]{DaLaMu21}). Next, we find convenient to introduce the notation
\[
M^{\sharp}_{jl}[f](x)\equiv \tilde{\nu}_{j}(x)\frac{\partial f}{\partial x_{l}}(x)
-
\tilde{\nu}_{l}(x)\frac{\partial f}{\partial x_{j}}(x)\qquad\forall x\in \overline{\Omega}\,,
\]
for all $f\in C^{1}(\overline{\Omega})$. If necessary, we write
$M^{\sharp}_{jl,x}$ to emphasize that we are taking $x$ as variable of the differential operator $M^{\sharp}_{jl}$. Next we fix $x\in \Omega$ and we note that
\begin{eqnarray}\label{thm:slay1}
\lefteqn{
M^{\sharp}_{jl}\left[v^+_\Omega[S_{ {\mathbf{a}} },\mu]\right](x)
}
\\ \nonumber
&&\qquad
=\tilde{\nu}_j(x)\frac{\partial}{\partial x_l}v^+_\Omega[S_{ {\mathbf{a}} },\mu](x)
-\tilde{\nu}_l(x)\frac{\partial}{\partial x_j}v^+_\Omega[S_{ {\mathbf{a}} },\mu](x)
\\ \nonumber
&&\qquad
= \int_{\partial\Omega}(\tilde{\nu}_j(x)-\nu_j(y))\frac{\partial}{\partial x_l}S_{{\mathbf{a}}}(x-y) \mu(y)
\\ \nonumber
&&\qquad\quad
-(\tilde{\nu}_l(x)-\nu_l(y))\frac{\partial}{\partial x_j}S_{{\mathbf{a}}}(x-y)  \mu(y)
\,d\sigma_y
\\ \nonumber
&&\qquad\quad
+ \int_{\partial\Omega}\left[\nu_j(y)\frac{\partial}{\partial x_l}\left(S_{{\mathbf{a}}}(x-y)\right) 
-\nu_l(y)\frac{\partial}{\partial x_j}\left(S_{{\mathbf{a}}}(x-y)\right)\right]\mu(y)
\,d\sigma_y
\\ \nonumber
&&\qquad
= \int_{\partial\Omega}(\tilde{\nu}_j(x)-\tilde{\nu}_j(y))\frac{\partial}{\partial x_l}S_{{\mathbf{a}}}(x-y) \mu(y)
\\ \nonumber
&&\qquad\quad
-(\tilde{\nu}_l(x)-\tilde{\nu}_l(y))\frac{\partial}{\partial x_j}S_{{\mathbf{a}}}(x-y)  \mu(y)
\,d\sigma_y
\\ \nonumber
&&\qquad\quad
 -\int_{\partial\Omega}
\left[\tilde{\nu}_j(y)\frac{\partial}{\partial y_l}S_{{\mathbf{a}}}(x-y)
-
\tilde{\nu}_l(y)
\frac{\partial}{\partial y_j}S_{{\mathbf{a}}}(x-y)
\right] \mu(y)
\,d\sigma_y
\\ \nonumber
&&\qquad
= \int_{\partial\Omega}(\tilde{\nu}_j(x)-\tilde{\nu}_j(y))\frac{\partial}{\partial x_l}S_{{\mathbf{a}}}(x-y) \mu(y)
\\ \nonumber
&&\qquad\quad
-(\tilde{\nu}_l(x)-\tilde{\nu}_l(y))\frac{\partial}{\partial x_j}S_{{\mathbf{a}}}(x-y)  \mu(y)
\,d\sigma_y
\\ \nonumber
&&\qquad\quad
-\int_{\partial\Omega}
  M_{jl,y}
\left[ S_{{\mathbf{a}}}(x-y)
\right] \mu(y)
\,d\sigma_y
\\ \nonumber
&&\qquad
= \int_{\partial\Omega}(\tilde{\nu}_j(x)-\tilde{\nu}_j(y))\frac{\partial}{\partial x_l}S_{{\mathbf{a}}}(x-y) \mu(y)
\\ \nonumber
&&\qquad\quad
-(\tilde{\nu}_l(x)-\tilde{\nu}_l(y))\frac{\partial}{\partial x_j}S_{{\mathbf{a}}}(x-y)  \mu(y)
\,d\sigma_y
\\ \nonumber
&&\qquad\quad
+ \int_{\partial\Omega}S_{{\mathbf{a}}}(x-y) 
M_{jl}[\mu](y)
\,d\sigma_y 
\end{eqnarray}
(cf. Lemma \ref{lem:gagrelip}). Since $(\tilde{\nu},\mu)\in C^{0,1}(\overline{\Omega},{\mathbb{R}}^n)\times L^\infty (\partial\Omega)$ the first integral in the   right hand side of (\ref{thm:slay1}) defines a continuous function of $x\in\overline{\Omega}$ (cf.~\cite[Thm. 8.1 (i)]{DoLa17}). Since $\Omega$ is of class $C^1$ and $M_{lj}[\mu]\in L^\infty(\partial\Omega)$, the second integral in the right hand side of (\ref{thm:slay1}) defines a continuous function of $x\in \overline{\Omega}$
 (cf.~\cite[Lem.~4.2 (i) with $G=\partial\Omega$, Lem.~6.2]{DoLa17}). Since $v^+_\Omega[S_{ {\mathbf{a}} },\mu]$ is of class $C^1(\overline{\Omega})$, then $M^{\sharp}_{lj}\left[v^+_\Omega[S_{ {\mathbf{a}} },\mu]\right]$ is continuous on $\overline{\Omega}$. Hence the left and right hand sides of (\ref{thm:slay1}) must be equal for all $x\in \overline{\Omega}$ and thus we have
 \begin{eqnarray}\label{thm:slay2}
\lefteqn{
M_{jl}[V_\Omega[S_{{\mathbf{a}}} ,\mu]](x)=M^{\sharp}_{jl}\left[v^+_\Omega[S_{ {\mathbf{a}} },\mu]\right](x)
}
\\ \nonumber
&&\qquad
= \int_{\partial\Omega}(\nu_j(x)-\nu_j(y))\frac{\partial}{\partial x_l}S_{{\mathbf{a}}}(x-y) \mu(y)
\\ \nonumber
&&\qquad\quad
-( \nu_l(x)-\nu_l(y))\frac{\partial}{\partial x_j}S_{{\mathbf{a}}}(x-y)  \mu(y)
\,d\sigma_y
\\ \nonumber
&&\qquad\quad
+ \int_{\partial\Omega}S_{{\mathbf{a}}}(x-y) 
M_{jl}[\mu](y)
\,d\sigma_y
\\ \nonumber
&&\qquad
=Q_{l}[\nu_j,\mu](x)-Q_{j}[\nu_l,\mu](x)+V_\Omega[S_{{\mathbf{a}}} ,M_{jl}[\mu]](x)
\qquad\forall x\in \partial\Omega\,.
\end{eqnarray}
Since the components of $\nu$ are of class $C^{0,1}$,    the first two terms in the right hand side of (\ref{thm:slay2}) define linear and continuous maps of the variable $\mu$ from $L^\infty(\partial\Omega)$ to $C^{0,\omega_1(\cdot)}(\partial\Omega)$
  (cf.~\cite[Thm.~8.2 (i)]{DoLa17}).
Since $M_{lj}$ is continuous from $C^{0,1}(\partial\Omega)$ to   $L^\infty(\partial\Omega)$ and $V_\Omega[S_{{\mathbf{a}}} ,\cdot]$ is linear and continuous from $L^\infty(\partial\Omega)$ to
 $C^{0,\omega_1(\cdot)}(\partial\Omega)$ (cf.~\cite[Thm.~7.2]{DoLa17}, \cite[Prop.~5.3]{La22d}), we conclude that the right hand side of (\ref{thm:slay2}) defines a linear and continuous map of the variable $\mu$ from $C^{0,1}(\partial\Omega)$ to $C^{0,\omega_1(\cdot)}(\partial\Omega)$. Hence equality (\ref{thm:slay2}) implies the validity of statement (jjj) and the proof is complete.
 
 Next we assume that the statement holds for $m\geq 1$, and we prove it for $m+1$. By the definition of norm in $C^{m+1,\omega_1(\cdot)}(\partial\Omega)$, it suffices to show that
\begin{enumerate}
\item[(l)] $V_\Omega[S_{{\mathbf{a}}} ,\cdot]$ is linear and continuous from  $C^{m,1}(\partial\Omega)$ to $C^{0}(\partial\Omega)$.
\item[(ll)] $V_\Omega[S_{{\mathbf{a}}} ,\mu]$ is continuously differentiable on $\partial\Omega$ for all $\mu$ in $C^{m,1}(\partial\Omega)$.
\item[(lll)] $M_{jl}[V_\Omega[S_{{\mathbf{a}}} ,\cdot]]$ is linear and continuous from  $C^{m,1}(\partial\Omega)$ to $C^{m,\omega_1(\cdot)}(\partial\Omega)$ for all $j,l\in \{1,\dots,n\}$,
 \end{enumerate}
 (cf.~\textit{e.g.}, \cite[Lem.~2.3]{DoLa17}). Since $C^{m,1}(\partial\Omega)$ is continuously imbedded into $C^{0,1}(\partial\Omega)$, statements (l), (ll) hold by case $m=1$.  We now prove statement (lll) by exploiting formula (\ref{thm:slay2}) and the inductive assumption. Since the components of $\nu$ are of class $C^{m,1}$ and $\Omega$ is of class $C^{m+1,1}$ and accordingly of class $C^{m+1,\alpha}$ for all $\alpha\in]0,1[$,    the first two terms in the right hand side of (\ref{thm:slay2}) define linear and continuous maps of the variable $\mu$ from $C^{m}(\partial\Omega)$ to $C^{m,\omega_1(\cdot)}$
  (cf.~\cite[Thm.~8.3 (i)]{DoLa17}).
  
  Since $M_{jl}$ is continuous from $C^{m,1}(\partial\Omega)$ to $C^{m-1,1}(\partial\Omega)$ and  the  inductive assumption implies that
   $V_\Omega[S_{{\mathbf{a}}} ,\cdot]$ is linear and continuous from $C^{m-1,1}(\partial\Omega)$ to
 $C^{m,\omega_1(\cdot)}$, we conclude that the right hand side of (\ref{thm:slay2}) defines a linear and continuous map of the variable $\mu$ from $C^{m,1}(\partial\Omega)$ to $C^{m,\omega_1(\cdot)}(\partial\Omega)$. Hence equality (\ref{thm:slay2}) implies the validity of statement (lll) and the proof is complete.\hfill  $\Box$ 

\vspace{\baselineskip}

\section{Analysis of the map $Q_r$}

We are now ready to prove the following statement.

\begin{thm}
\label{thm:qrsm}
Let ${\mathbf{a}}$ be as in (\ref{introd0}), (\ref{ellip}), (\ref{symr}). Let $S_{ {\mathbf{a}} }$ be a fundamental solution of $P[{\mathbf{a}},D]$. 
   Let $\alpha\in]0,1]$. 
Let $\Omega$ be a bounded open subset of ${\mathbb{R}}^{n}$ of class $C^{2,\alpha}$. Let $j\in\{1,\dots,n\}$. Then   the following statements hold.
\begin{enumerate}
\item[(i)]  If $\alpha\in]0,1[$, $\beta\in]0,\alpha]$, then
the  bilinear map $Q_j\left[\cdot,\cdot\right]$  from  the space $C^{1,\alpha}(\partial\Omega)\times C^{0,\beta}(\partial\Omega)$ to $C^{1,\beta}(\partial\Omega)$ which takes a pair   $(g,\mu)$ to $Q_j\left[ g,\mu\right]$ is continuous (cf.~(\ref{qrs0})). 
\item[(ii)] If $\alpha=1$, $\beta=1$, then
the  bilinear map $Q_j\left[\cdot,\cdot\right]$  from  the space $C^{1,1}(\partial\Omega)\times C^{0,1}(\partial\Omega)$ to $C^{1,\omega_1(\cdot)}(\partial\Omega)$ which takes a pair   $(g,\mu)$ to $Q_j\left[ g,\mu\right]$ is continuous (cf.~(\ref{qrs0})). 
\end{enumerate}
\end{thm}
{\bf Proof.} We prove statements (i) and (ii) at the same time and make some appropriate comment when the two proofs present some difference. By the definition of norm in $C^{1,\beta}(\partial\Omega)$ with $\beta\in]0,1[$ and in $C^{1,\omega_1(\cdot)}(\partial\Omega)$, it suffices to show that
\begin{enumerate}
\item[(j)] $Q_j$ is bilinear and continuous from  $C^{1,\alpha}(\partial\Omega)\times C^{0,\beta}(\partial\Omega)$ to $C^{0}(\partial\Omega)$.
\item[(jj)] $Q_j\left[ g,\mu\right]$ is continuously differentiable on $\partial\Omega$ for all $(g,\mu) $ in $C^{1,\alpha}(\partial\Omega)\times C^{0,\beta}(\partial\Omega)$.
\item[(jjj)] ${\mathrm{grad}}_{\partial\Omega}Q_j\left[\cdot,\cdot\right]$ is bilinear and continuous from  $C^{1,\alpha}(\partial\Omega)\times C^{0,\beta}(\partial\Omega)$ to $C^{0,\beta}(\partial\Omega)$ in case of statement (i) and to $C^{0,\omega_1(\cdot)}(\partial\Omega)$ in case of statement (ii), 
\end{enumerate}
(cf.~\textit{e.g.}, \cite[Lem.~2.3]{DoLa17}). For a proof of (j), we refer to \cite[Thm.~8.2 (i)]{DoLa17}. Next we set
\[
K_j[g](x,y)\equiv (g(x)-g(y))\frac{\partial S_{ {\mathbf{a}} }}{\partial x_{j}}(x-y)
\]
for all $(x,y)\in (\partial\Omega)^2\setminus {\mathbb{D}}_{\partial\Omega}$  and $g\in C^{1,\alpha}(\partial\Omega)$ and we note that
\[
Q_j\left[ g,\mu\right](x)=\int_{\partial\Omega}K_j(g)(x,y)\mu(y)\,d\sigma_y
\qquad\forall (g,\mu)\in C^{1,\alpha}(\partial\Omega)\times C^{0,\beta}(\partial\Omega) 
\]
and we turn to the proof of (jj). To do so, we resort to a classical differentiation Theorem in the form of \cite[Thm.~6.2]{La22b} and we turn to verify its assumptions. It is known that
\[
\frac{\partial S_{ {\mathbf{a}} }}{\partial x_{j}}(x-y)\in {\mathcal{K}}_{n-1,n,1}((\partial\Omega)\times(\partial\Omega))\subseteq  {\mathcal{K}}_{n-1,n-1+\alpha,\alpha}((\partial\Omega)\times(\partial\Omega))
 \]
(cf.~Lemma \ref{lem:kelem}, \cite[Lem.~4.3]{DoLa17}). Then the product Lemma \cite[Lem.~3.4 (ii)]{La22b} implies that 
\[
K_j[g]\in {\mathcal{K}}_{n-1-1,n-1,1-(1-1)} ((\partial\Omega)\times(\partial\Omega))={\mathcal{K}}_{n-2,n-1,1} ((\partial\Omega)\times(\partial\Omega))
\]
for all  $g\in C^{0,1}(\partial\Omega)$ and that there exists $c_1\in]0,+\infty[$ such that
\begin{equation}\label{thm:qrsm1}
\|K_j[g]\|_{ {\mathcal{K}}_{n-2,n-1,1}  ((\partial\Omega)\times(\partial\Omega))}\leq c_1\left\|\frac{\partial S_{ {\mathbf{a}} }}{\partial x_{j}}(x-y)\right\|_{ {\mathcal{K}}_{n-1,n,1} ((\partial\Omega)\times(\partial\Omega))}\|g\|_{C^{0,1}(\partial\Omega)}
\end{equation}
for all $g\in C^{0,1}(\partial\Omega)$. Since $\Omega$ is of class $C^{2,\alpha}$ and $(g,1)$ belongs to $C^{1,\alpha}(\partial\Omega)\times C^{1,\beta}(\partial\Omega) $, Theorem 8.3 of \cite{DoLa17} implies that
\begin{equation}\label{thm:qrsm2}
\int_{\partial\Omega}K_j[g](\cdot,y)\,d\sigma_y=Q_j[g,1]\in
\left\{
\begin{array}{ll}
 C^{1,\alpha}(\partial\Omega)& \text{if}\ \alpha\in]0,1[\,,
 \\
 C^{1,\omega_1(\cdot)}(\partial\Omega)& \text{if}\ \alpha=1\,,
\end{array}
\right.  
\end{equation}
for all $g\in C^{1,\alpha}(\partial\Omega)$ and that  
\begin{eqnarray}\label{thm:qrsm2a}
\lefteqn{
Q_j[\cdot,1]\ \text{is\ linear\ and\ continuous\ from }\ C^{1,\alpha}(\partial\Omega)
}
\\ \nonumber
&&\qquad\qquad\qquad\qquad\qquad
  \text{to}\ 
\left\{
\begin{array}{ll}
 C^{1,\alpha}(\partial\Omega)& \text{if}\ \alpha\in]0,1[\,,
 \\
 C^{1,\omega_1(\cdot)}(\partial\Omega)& \text{if}\ \alpha=1\,.
\end{array}
\right.  
\end{eqnarray}
We also note that
\[
K_j[g]\in C^1((\partial\Omega)\setminus\{y\})\qquad\forall y\in\partial\Omega\,,
\]
for all $g\in C^{1,\alpha}(\partial\Omega)$.  Next we compute the tangential gradient with respect to $x$  of $K_j[g]$. By the Leibnitz rule, we have
\begin{eqnarray}\label{thm:qrsm3}
\lefteqn{
{\mathrm{grad}}_{\partial\Omega,x} K_j[g](x,y) 
}
\\ \nonumber
&&\qquad
={\mathrm{grad}}_{\partial\Omega,x} g(x)\frac{\partial S_{ {\mathbf{a}} }}{\partial x_{j}}(x-y)
+(g(x)-g(y)){\mathrm{grad}}_{\partial\Omega,x}\left(\frac{\partial S_{ {\mathbf{a}} }}{\partial x_{j}}(x-y)\right)
\end{eqnarray}
for all $(x,y)\in (\partial\Omega)^2\setminus {\mathbb{D}}_{\partial\Omega}$  and  $g\in C^{1,\alpha}(\partial\Omega)$. Since
$\frac{\partial S_{ {\mathbf{a}} }}{\partial x_{j}}(x-y)$ belongs to $  
 {\mathcal{K}}_{n-1,n-1+\alpha,\alpha}((\partial\Omega)\times(\partial\Omega))$, $n-1<n-1+\alpha$, and the components of 
${\mathrm{grad}}_{\partial\Omega,x} g$ are $\alpha$-H\"{o}lder continuous, the product Proposition \ref{prop:prkerho} (ii) implies that
\begin{equation}\label{thm:qrsm4}
({\mathrm{grad}}_{\partial\Omega,x} g)_h(x)\frac{\partial S_{ {\mathbf{a}} }}{\partial x_{j}}(x-y)\in {\mathcal{K}}_{n-1,n-1+\alpha,\alpha}((\partial\Omega)\times(\partial\Omega)) 
\end{equation}
for all $g\in C^{1,\alpha}(\partial\Omega)$ and $h\in\{1,\dots,n\}$.  Since $g$ is Lipschitz continuous, Lemma \ref{lem:hoker} imples that 
\begin{equation}\label{thm:qrsm5}
(g(x)-g(y))\in {\mathcal{K}}_{-1,0,1}((\partial\Omega)\times(\partial\Omega)) 
\end{equation} 
for all $g\in C^{1,\alpha}(\partial\Omega)$. Since $\Omega$ is of class $C^{2,\alpha}$, then it is also of class $C^{1,1}$ and  Lemma \ref{lem:tgdjsgen} implies that
\begin{equation}\label{thm:qrsm6}
\left({\mathrm{grad}}_{\partial\Omega,x} \left(\frac{\partial S_{{\mathbf{a}}}}{\partial x_j}(x-y) \right)\right)_h\in {\mathcal{K}}_{n,n+1,1}((\partial\Omega)\times (\partial\Omega)) 
\end{equation}
for all $h\in\{1,\dots,n\}$. Then the product Theorem \ref{thm:kerpro} (ii) implies that
\begin{equation}\label{thm:qrsm7}
(g(x)-g(y))\left({\mathrm{grad}}_{\partial\Omega,x} \left(\frac{\partial S_{{\mathbf{a}}}}{\partial x_j}(x-y)\right) \right)_h\in  {\mathcal{K}}_{n-1,n,1}((\partial\Omega)\times (\partial\Omega)) 
\end{equation}
and that there exists $c_2\in]0,+\infty[$ such that
\begin{eqnarray}\label{thm:qrsm7a}
\lefteqn{
\left\|(g(x)-g(y))\left({\mathrm{grad}}_{\partial\Omega,x} \left(\frac{\partial S_{{\mathbf{a}}}}{\partial x_j}(x-y)\right) \right)_h\right\|_{{\mathcal{K}}_{n-1,n,1}((\partial\Omega)\times (\partial\Omega)) }
}
\\ \nonumber
&&\qquad\qquad\qquad\qquad\qquad\qquad\qquad
\leq c_2\|g\|_{C^{0,1}(\partial\Omega)}\qquad\forall g\in C^{0,1}(\partial\Omega)
\end{eqnarray}
for all   $h\in\{1,\dots,n\}$.
In particular, equality (\ref{thm:qrsm3}) and the memberships of (\ref{thm:qrsm4}), (\ref{thm:qrsm7}) imply that 
\[
({\mathrm{grad}}_{\partial\Omega,x} K_j[g])_h\in {\mathcal{K}}_{n-1,(\partial\Omega)\times(\partial\Omega)}\qquad\forall g\in C^{1,\alpha}(\partial\Omega)  
\]
for all $h\in\{1,\dots,n\}$. Then $\int_{\partial\Omega}K_j[g](\cdot,y)\mu(y)\,d\sigma_y$ is continuously differentiable and
\begin{eqnarray}\label{thm:qrsm8}
\lefteqn{
 {\mathrm{grad}}_{\partial\Omega} \int_{\partial\Omega}K_j[g](x,y)\mu(y)\,d\sigma_y
}
\\ \nonumber
&&\qquad
=\int_{\partial\Omega}[{\mathrm{grad}}_{\partial\Omega,x}K_j[g](x,y)](\mu(y)-\mu(x))\,d\sigma_y
\\ \nonumber
&&\qquad\quad
+\mu(x){\mathrm{grad}}_{\partial\Omega} \int_{\partial\Omega}K_j[g](x,y) \,d\sigma_y \,,
\end{eqnarray}
for all $x\in \partial\Omega$ and for all $(g,\mu)\in C^{1,\alpha}(\partial\Omega)\times C^{0,\beta}(\partial\Omega)$ (cf.~\cite[Thm.~6.2]{La22b}) and the proof of (jj) is complete. We now turn to prove (jjj). By equalities (\ref{thm:qrsm3}) and (\ref{thm:qrsm8}), we have
\begin{eqnarray}\label{thm:qrsm9}
\lefteqn{
({\mathrm{grad}}_{\partial\Omega} \int_{\partial\Omega}K_j[g](x,y)\mu(y)\,d\sigma_y)_h
}
\\ \nonumber
&&\quad
= 
({\mathrm{grad}}_{\partial\Omega,x} g(x))_h\int_{\partial\Omega}(\mu(y)-\mu(x))\frac{\partial S_{ {\mathbf{a}} }}{\partial x_{j}}(x-y) \,d\sigma_y 
\\ \nonumber
&&\quad\quad
+\int_{\partial\Omega}(g(x)-g(y))\left({\mathrm{grad}}_{\partial\Omega,x}\left(\frac{\partial S_{ {\mathbf{a}} }}{\partial x_{j}}(x-y)\right)\right)_h(\mu(y)-\mu(x)) \,d\sigma_y
\\ \nonumber
&&\quad\quad
+\mu(x)({\mathrm{grad}}_{\partial\Omega}  \int_{\partial\Omega}K_j[g](x,y) \,d\sigma_y)_h \,,
\end{eqnarray}
for all $x\in \partial\Omega$, for all $(g,\mu)\in C^{1,\alpha}(\partial\Omega)\times C^{0,\beta}(\partial\Omega)$ and $h\in \{1,\dots,n\}$. In order to prove statement (jjj) it suffices to show that each addendum in the right hand side of formula (\ref{thm:qrsm9}), defines a bilinear and continuous map from  $C^{1,\alpha}(\partial\Omega)\times C^{0,\beta}(\partial\Omega)$ to $C^{0,\beta}(\partial\Omega)$ in case of statement (i) and to $C^{0,\omega_1(\cdot)}(\partial\Omega)$ in case of statement (ii). We first consider the first addendum. 
Since $\Omega$ is of class $C^{2,\alpha}$, \cite[Thm.~8.2]{DoLa17}
 implies that $Q_j[\cdot,1]$ is linear and continuous from
$C^{0,\beta}(\partial\Omega)$ to $C^{0,\beta}(\partial\Omega)$ in case of  statement (i)  and from $C^{0,1}(\partial\Omega)$ to $C^{0,\omega_1(\cdot)}(\partial\Omega)$ in case of  statement (ii).  
 Since the components of ${\mathrm{grad}}_{\partial\Omega,x}$ are linear and continuous from $C^{1,\alpha}(\partial\Omega)$ to $C^{0,\alpha}(\partial\Omega)$
and the pointwise product is bilinear and continuous in (generalized) H\"{o}lder spaces (cf.~\textit{e.g.}, \cite[Lem.~2.5]{DoLa17}), we deduce that the first addendum in the right hand side of formula (\ref{thm:qrsm9}), defines a bilinear and continuous map from  $C^{1,\alpha}(\partial\Omega)\times C^{0,\beta}(\partial\Omega)$ to $C^{0,\beta}(\partial\Omega)$ in case of  statement (i) and to $C^{0,\omega_1(\cdot)}(\partial\Omega)$ in case of  statement (ii).

Next we consider the second addendum  in the right hand side of formula (\ref{thm:qrsm9}), that is an integral operator with the kernel of 
(\ref{thm:qrsm7}). We plan to apply a result of \cite[Prop.~6.3 (ii)]{La22a}. Since $Y\equiv\partial\Omega$ is a compact manifold of class $C^1$ that is imbedded in ${\mathbb{R}}^n$,  $Y$ can be proved to be  strongly upper $(n-1)$-Ahlfors regular with 
respect to $Y$ in the sense of   \cite[(1.4)]{La22a}. %(see \cite{La22e}). 
Then we set
\[
s_1\equiv n-1\,,\quad
s_2\equiv n\,,\quad
s_3 \equiv  1\,.
\]
and we note that
\[
(n-1)>0 \,,\ \ \beta\in]0,\alpha]\subseteq]0,1]\,,\ \ 
s_1\in [\beta,(n-1)+\beta[\,,\ \  s_2\in [\beta,+\infty[ 
\]
and that
\begin{eqnarray*}
\lefteqn{s_2-\beta=n-\beta>n-1\,,}
\\ \nonumber
&&\qquad\qquad\qquad
 s_2=n<n-1+\beta+1=n-1+\beta+s_3
\quad \text{if}\ \beta<1\,,
\\
\lefteqn{s_2-\beta=n-\beta=n-1 \quad \text{if}\ \beta=1\,.}
\end{eqnarray*}
Then \cite[Prop.~6.3 (ii) (b) and (bb)]{La22a} implies that the map 
\[
\text{from}\ {\mathcal{K}}^\sharp_{n-1,n,1}((\partial\Omega)\times (\partial\Omega))\times C^{0,\beta}(\partial\Omega)\ 
\text{to}\ 
\left\{
\begin{array}{ll}
 C^{0,\beta}(\partial\Omega)& \text{if}\ \beta\in]0,1[\,,
 \\
 C^{1,\omega_1(\cdot)}(\partial\Omega)& \text{if}\ \beta=1\,.
\end{array}
\right.  \]
that takes a pair $(K,\mu)$ to $\int_{\partial\Omega}K(\cdot,y)(\mu(y)-\mu(x))\,d\sigma_y$ is bilinear and continuous. Thus it suffices to show that the map  
\[
\text{from}\ C^{1,\alpha}(\partial\Omega)\ \text{to}\  {\mathcal{K}}^\sharp_{n-1,n,1}((\partial\Omega)\times (\partial\Omega))
\]
 that takes $g$ to the kernel in (\ref{thm:qrsm7}) is linear and continuous. By (\ref{thm:qrsm7a}) we know that
such a map is linear and continuous from $C^{1,\alpha}(\partial\Omega)$ to the space ${\mathcal{K}}_{n-1,n,1}((\partial\Omega)\times (\partial\Omega))$. Then by Lemma \ref{lem:lipgw} of the Appendix, there exists $c_{\Omega,1}\in]0,+\infty[$ such that
\begin{eqnarray} \nonumber
\lefteqn{
\sup_{x\in \partial\Omega}\sup_{r\in ]0,+\infty[}
\left|
\int_{(\partial\Omega)\setminus {\mathbb{B}}_n(x,r)}
(g(x)-g(y))\left({\mathrm{grad}}_{\partial\Omega,x}\left(\frac{\partial S_{ {\mathbf{a}} }}{\partial x_{j}}(x-y)\right)\right)_h
 \,d\sigma_y
\right|
}
\\  \nonumber
&&
\leq\sup_{x\in \partial\Omega}\sup_{r\in ]0,+\infty[}
\biggl|
\int_{(\partial\Omega)\setminus {\mathbb{B}}_n(x,r)}
(g(x)-g(y)
+
{\mathrm{grad}}_{\partial\Omega}g(x)(y-x)
)
\\  \label{thm:qrsm10}
&&\quad
\times
\left({\mathrm{grad}}_{\partial\Omega,x}\left(\frac{\partial S_{ {\mathbf{a}} }}{\partial x_{j}}(x-y)\right)\right)_h
 \,d\sigma_y
\biggr|
\\  \nonumber
&&
+\sup_{x\in \partial\Omega}\sup_{r\in ]0,+\infty[}
\biggl|
\int_{(\partial\Omega)\setminus {\mathbb{B}}_n(x,r)}
({\mathrm{grad}}_{\partial\Omega}g(x)(y-x)
)
\\  \nonumber
&&\quad
\times
({\mathrm{grad}}_{\partial\Omega,x}\frac{\partial S_{ {\mathbf{a}} }}{\partial x_{j}}(x-y))_h
 \,d\sigma_y
\biggr|
\\  \nonumber
&&
\leq 
c_{\Omega,1}
\biggl(\sup_{\partial\Omega}|g|
+\sup_{\partial\Omega}|{\mathrm{grad}}_{\partial\Omega}g|
\\  \nonumber
&&\quad
+|{\mathrm{grad}}_{\partial\Omega}g:\partial\Omega|_\alpha\biggr)
\left\|
\left({\mathrm{grad}}_{\partial\Omega,x}\left(\frac{\partial S_{ {\mathbf{a}} }}{\partial x_{j}}(x-y)\right)\right)_h
\right\|_{{\mathcal{K}}_{n,n+1,n}((\partial\Omega)\times (\partial\Omega))}
\\  \nonumber
&&\quad\qquad\qquad\times
\sup_{x\in \partial\Omega} 
 \int_{\partial\Omega}|x-y|^{1+\alpha-n}\,d\sigma_y
\\  \nonumber
&&\quad
+\sum_{z=1}^n\sup_{x\in \partial\Omega}|({\mathrm{grad}}_{\partial\Omega}g(x))_z|
\\  \nonumber
&&\quad
\times 2n \sup_{s\in\{1,\dots,n\}}\sup_{x\in \partial\Omega}
\sup_{r\in ]0,+\infty[}
\left|\,
\int_{(\partial\Omega)\setminus{\mathbb{B}}_{n}(x,r)}
(x_z-y_z)\frac{\partial}{\partial x_s}\frac{\partial S_{{\mathbf{a}}}}{\partial x_j}(x-y) 
\,d\sigma_y\right|
\end{eqnarray}
for all $g\in C^{1,\alpha}(\partial\Omega)$ and
$h\in\{1,\dots,n\}$ (see also (\ref{thm:maggnfs2}), (\ref{thm:qrsm6}) and Lemma \ref{lem:mahgnfs}). 
Hence, we deduce that the second addendum in the right hand side of formula (\ref{thm:qrsm9}), defines a bilinear and continuous map from  $C^{1,\alpha}(\partial\Omega)\times C^{0,\beta}(\partial\Omega)$ to $C^{0,\beta}(\partial\Omega)$ in case of statement (i) and to $C^{0,\omega_1(\cdot)}(\partial\Omega)$ to   in case of statement (ii).

Next we consider the third addendum in the right hand side of formula (\ref{thm:qrsm9}). Since  ${\mathrm{grad}}_{\partial\Omega,x}  $ is linear and continuous from $C^{1,\alpha}(\partial\Omega)$ to $C^{0,\alpha}(\partial\Omega)$, the continuity of $Q_j[\cdot,1]$ as in (\ref{thm:qrsm2a}) and the continuity of the pointwise product in generalized H\"{o}lder spaces (cf.~\textit{e.g.}, \cite[Lem.~2.5]{DoLa17})  imply that the third addendum in the right hand side of formula (\ref{thm:qrsm9}), defines a bilinear and continuous map from  $C^{1,\alpha}(\partial\Omega)\times C^{0,\beta}(\partial\Omega)$ to $C^{0,\beta}(\partial\Omega)$ in case of statement (i) and to $C^{0,\omega_1(\cdot)}(\partial\Omega)$ to   in case of statement (ii).
Hence, the proof of (jjj) and of the theorem is complete.\hfill  $\Box$ 

\vspace{\baselineskip}

In the previous theorem, we have considered sets of class $C^{2,\alpha}$. We are now ready to consider case $C^{m,\alpha}$ by an inductive argument on $m$ as in the proof of   \cite[Thm.~8.3]{DoLa17}.
\begin{thm}
\label{thm:qrsmm}
Let ${\mathbf{a}}$ be as in (\ref{introd0}), (\ref{ellip}), (\ref{symr}). Let $S_{ {\mathbf{a}} }$ be a fundamental solution of $P[{\mathbf{a}},D]$. 
 Let $\alpha\in]0,1]$.  Let $m\in {\mathbb{N}}$, $m\geq 2$.
Let $\Omega$ be a bounded open subset of ${\mathbb{R}}^{n}$ of class $C^{m,\alpha}$. Let $r\in\{1,\dots,n\}$. 
Then the following statements hold.
 \begin{enumerate}
\item[(i)]  If $\alpha\in]0,1[$ and $\beta\in]0,\alpha]$, then
the  bilinear map $Q_r\left[\cdot,\cdot\right]$  from  the space $C^{m-1,\alpha}(\partial\Omega)\times C^{m-2,\beta}(\partial\Omega)$ to $C^{m-1,\beta}(\partial\Omega)$  which takes a pair   $(g,\mu)$ to $Q_r\left[ g,\mu\right]$ is continuous (cf.~(\ref{qrs0})). 
\item[(ii)]  If $\alpha=1$ and $\beta=1$, then
the  bilinear map $Q_r\left[\cdot,\cdot\right]$  from  the space $C^{m-1,1}(\partial\Omega)\times C^{m-2,1}(\partial\Omega)$ to $C^{m-1,\omega_1(\cdot)}(\partial\Omega)$ which takes a pair   $(g,\mu)$ to $Q_r\left[ g,\mu\right]$ is continuous (cf.~(\ref{qrs0})). 
\end{enumerate}
\end{thm}
{\bf Proof.} We prove statements (i) and (ii) at the same time and make some appropriate comment when the two proofs present some difference. We proceed by induction on $m$. Case $m=2$ holds by Theorem \ref{thm:qrsm}. We now prove that if the statement holds for $m$, then it holds also for $m+1$. Then we now assume that $\Omega$ is of class $C^{m+1,\alpha}$ and we prove that 
  $Q_r\left[\cdot,\cdot\right]$ is continuous from  $C^{m,\alpha}(\partial\Omega)\times C^{m-1,\beta}(\partial\Omega)$ to $C^{m,\beta}(\partial\Omega)$  in case of  statement (i) and to $C^{m,\omega_1(\cdot)}(\partial\Omega)$ in case of  statement (ii).\par
  
By the definition of norm in $C^{m,\beta}(\partial\Omega)$ with $\beta\in]0,1[$ and in 
  $C^{m,\omega_1(\cdot)}(\partial\Omega)$, it suffices to show that
\begin{enumerate}
\item[(j)] $Q_r$ is bilinear and continuous from  $C^{m,\alpha}(\partial\Omega)\times C^{m-1,\beta}(\partial\Omega)$ to $C^{0}(\partial\Omega)$.
\item[(jj)] $Q_r\left[ g,\mu\right]$ is continuously differentiable for all $(g,\mu)$ in $C^{m,\alpha}(\partial\Omega)\times C^{m-1,\beta}(\partial\Omega)$.
\item[(jjj)] If  $j,l\in \{1,\dots,n\}$, then $M_{lj}\left[Q_r\left[\cdot,\cdot\right]\right]$ is bilinear and continuous from  $C^{m,\alpha}(\partial\Omega)\times C^{m-1,\beta}(\partial\Omega)$ to $C^{m-1,\beta}(\partial\Omega)$ in case of  statement (i) and to $C^{m-1,\omega_1(\cdot)}(\partial\Omega)$ in case of  statement (ii),
\end{enumerate}
(cf.~\textit{e.g.}, \cite[Lem.~2.3]{DoLa17}). Statements (j), (jj) hold  by the continuous embedding of $C^{m,\alpha}(\partial\Omega)\times C^{m-1,\beta}(\partial\Omega)$ into $C^{1,\alpha}(\partial\Omega)\times C^{0,\beta}(\partial\Omega)$ and by case $m=2$. We now prove statement (jjj).
 We first note that if $(g,\mu)$ belongs to $C^{m,\alpha}(\partial\Omega)\times C^{m-1,\beta}(\partial\Omega)$, then assumption $m\geq 2$ and \cite[Lem.~8.1]{DoLa17} imply that the following formula holds
\[
M_{lj}\left[Q_r\left[g,\mu \right]\right]=
P_{ljr}[g,\mu]\,,
\]
where
\begin{eqnarray*}
\lefteqn{
P_{ljr}[g,\mu](x)
}
\\ \nonumber
&& 
\equiv \biggl\{\nu_l(x)Q_r\left[ ({\mathrm{grad}}_{\partial\Omega}g)_j ,\mu\right](x) - \nu_j(x)Q_r\left[({\mathrm{grad}}_{\partial\Omega}g)_l ,\mu\right](x) \biggr\}
\\ \nonumber
&& 
 + \biggl\{\nu_l(x)Q_r\left[  g , 
 \sum_{s=1}^{n} M_{sj}[\sum_{h=1}^{n} \frac{ a_{sh}\nu_h}{
 \nu^{t}a^{(2)}\nu  }\mu]\right](x)
 \\ \nonumber
&& \qquad
- 
\nu_j(x)Q_r\left[ g , \sum_{s =1}^{n}M_{sl}[\sum_{h=1}^{n} \frac{a_{sh}\nu_h}{
  \nu^{t}a^{(2)}\nu} \mu ]\right](x)\biggr\}
\\ \nonumber
&&   
  + \sum_{s,h=1}^{n} a_{sh} \nu_l(x) \biggl\{\biggr.
  Q_s\left[\nu_j,\frac{M_{hr}[g]\mu}{ \nu^{t}a^{(2)}\nu}\right](x) 
  \\ \nonumber
&& \qquad
+ 
Q_s\left[ g, M_{hr}[\frac{\nu_j\mu}
{\nu^{t}a^{(2)}\nu} ]\right](x) 
\biggl.\biggr\}
\\ \nonumber
&&   
- \sum_{s,h=1}^{n} a_{sh} \nu_j(x) \biggl\{\biggr.
Q_s\left[\nu_l,\frac{M_{hr}[g]\mu}{
\nu^{t}a^{(2)}\nu
}\right](x) 
\\ \nonumber
&& \qquad
+ Q_s\left[g, M_{hr}[\frac{\nu_l\mu}{ \nu^{t}a^{(2)}\nu} ]\right](x) \biggl.\biggr\} 
\\ \nonumber
&& 
 -\sum_{t=1}^{n} a_{s}\biggl\{\biggr.\nu_l(x)
 Q_s\left[g,\frac{\nu_j\nu_r}{\nu^{t}a^{(2)}\nu}\mu\right](x) 
 \\ \nonumber
&& \qquad
 -\nu_j(x)
 Q_s\left[ g,\frac{\nu_l\nu_r}{
\nu^{t}a^{(2)}\nu
}\mu\right](x) \biggl.\biggr\} 
\\ \nonumber
&& 
-a\left\{ 
g(x)\left[
\nu_l(x)v_\Omega[S_{ {\mathbf{a}} }  ,  
\frac{\nu_j\nu_r}{
\nu^{t}a^{(2)}\nu}\mu](x)
-
\nu_j(x)v_\Omega[S_{ {\mathbf{a}} } ,  
\frac{\nu_l\nu_r}{
\nu^{t}a^{(2)}\nu}\mu](x)
\right]\right.
\\ \nonumber
&& 
-
\left.
\left[
\nu_l(x)v_\Omega[S_{ {\mathbf{a}} } ,  
g\frac{\nu_j\nu_r}{
\nu^{t}a^{(2)}\nu}\mu](x)
-
\nu_j(x)v_\Omega[S_{ {\mathbf{a}} }   ,  
g\frac{\nu_l\nu_r}{
\nu^{t}a^{(2)}\nu}\mu](x)
\right]
\right\} \ \ \forall x\in\partial\Omega\,,
\end{eqnarray*}  
for all $(g,\mu)\in C^{m,\alpha}(\partial\Omega)\times C^{m-1,\beta}(\partial\Omega)$.   We first prove that if  $(g,\mu)$ belongs to $C^{m,\alpha}(\partial\Omega)\times C^{m-1,\beta}(\partial\Omega)$, then each term in the right hand side  of the equality that defines $P_{ljr}[g,\mu]$ belongs to  $C^{m-1,\beta}(\partial\Omega)$  in case of  statement (i) and to $C^{m-1,\omega_1(\cdot)}(\partial\Omega)$ in case of  statement (ii). Then the proof of the continuity of $P_{ljr}$ as in (jjj) follows the same lines and is accordingly omitted.

By the continuity of all  the components of  ${\mathrm{grad}}_{\partial\Omega}$   from $C^{m,\alpha}(\partial\Omega)$
to $C^{m-1,\alpha}(\partial\Omega)$ we have
$({\mathrm{grad}}_{\partial\Omega}g)_j \in C^{m-1,\alpha}(\partial\Omega)$.

By the continuity of the imbedding of $C^{m-1,\beta}(\partial\Omega)$ into $C^{m-2,\beta}(\partial\Omega)$ we have
$\mu\in C^{m-2,\beta}(\partial\Omega)$.

By by the inductive assumption on $Q_r$,  
$Q_r\left[ ({\mathrm{grad}}_{\partial\Omega}g)_j ,\mu\right]$ belongs to $
C^{m-1,\beta}(\partial\Omega)$ in case of  statement (i) and to $C^{m-1,\omega_1(\cdot)}(\partial\Omega)$ in case of  statement (ii).

By the membership of the components of $\nu$ in $C^{m,\alpha}(\partial\Omega)\subseteq
C^{m-1,\alpha}(\partial\Omega)
$, and by the continuity of the pointwise product
\[
 \text{from}\   
 C^{m-1,\alpha}(\partial\Omega)\times  C^{m-1,\beta}(\partial\Omega)\  \text{to}\  C^{m-1,\beta}(\partial\Omega)
 \]
in case of  statement (i) 
and 
\[
\text{from}\ 
 C^{m-1,\alpha}(\partial\Omega)\times  C^{m-1,\omega_1(\cdot)}(\partial\Omega)\ \text{to}\  C^{m-1,\omega_1(\cdot)}(\partial\Omega)
 \]
  in case of statement (ii) (cf.~\textit{e.g.}, \cite[Lems.~2.4, 2.5]{DoLa17}), the  sum in the first pair of braces in the right hand side  of the equality that defines $P_{ljr}[g,\mu]$ belongs to $C^{m-1,\beta}(\partial\Omega)$  in case of  statement (i) and to $C^{m-1,\omega_1(\cdot)}(\partial\Omega)$ in case of  statement (ii).  

For the remaining terms the argument is similar and thus we merely outline it. Since the components of $\nu $ belong to  $C^{m,\alpha}(\partial\Omega)\subseteq
C^{m-1,\alpha}(\partial\Omega)
$, the continuity of
the pointwise product in Schauder spaces implies that 
\[
\left(g , 
 \sum_{s=1}^{n} M_{sj}[\sum_{h=1}^{n} \frac{ a_{sh}\nu_h}{
 \nu^{t}a^{(2)}\nu  }\mu]\right)\in C^{m-1,\alpha}(\partial\Omega)\times C^{m-2,\beta}(\partial\Omega)\,.
 \]
  Then the inductive assumption on $Q_r $ ensures that 
  \[Q_r\left[  g , 
 \sum_{s=1}^{n} M_{sj}[\sum_{h=1}^{n} \frac{ a_{sh}\nu_h}{
 \nu^{t}a^{(2)}\nu  }\mu]\right]
 \]
 belongs to $C^{m-1,\beta}(\partial\Omega)$  in case of  statement (i) and to $C^{m-1,\omega_1(\cdot)}(\partial\Omega)$ in case of  statement (ii).  Then again the continuity of
the pointwise product in Schauder spaces implies that  
the sum in the second pair of braces in the right hand side  of the equality that defines
  $P_{ljr}[g,\mu]$  belongs to $C^{m-1,\beta}(\partial\Omega)$  in case of  statement (i) and to $C^{m-1,\omega_1(\cdot)}(\partial\Omega)$ in case of  statement (ii).\par

 By the membership of the components of $\nu$ in $C^{m,\alpha}(\partial\Omega)\subseteq
C^{m-1,\alpha}(\partial\Omega)
$,   by the continuity of $M_{hr}$ from $C^{m,\alpha}(\partial\Omega)$ to $C^{m-1,\alpha}(\partial\Omega)$ and by the continuity of the imbedding from $C^{m-1,\beta}(\partial\Omega)$ to $C^{m-2,\beta}(\partial\Omega)$,
by the continuity of
the pointwise product in Schauder spaces (cf.~\textit{e.g.}, \cite[Lems.~2.4, 2.5]{DoLa17}),  we have
\[
\left( \nu_j,\frac{M_{hr}[g]\mu}{ \nu^{t}a^{(2)}\nu}\right)
\in C^{m-1,\alpha}(\partial\Omega)\times C^{m-2,\beta}(\partial\Omega)\,.
 \]
Then the inductive assumption on $Q_s$ ensures that 
\[
Q_s\left[\nu_j,\frac{M_{hr}[g]\mu}{ \nu^{t}a^{(2)}\nu}\right]
\]
belongs to $C^{m-1,\beta}(\partial\Omega)$  in case of  statement (i) and to $C^{m-1,\omega_1(\cdot)}(\partial\Omega)$ in case of  statement (ii). Similarly, 
\[
Q_s\left[ g, M_{hr}[\frac{\nu_j\mu}
{\nu^{t}a^{(2)}\nu} ]\right]
\]
belongs to $C^{m-1,\beta}(\partial\Omega)$  in case of  statement (i) and to $C^{m-1,\omega_1(\cdot)}(\partial\Omega)$ in case of  statement (ii) and thus again by the continuity of the pointwise product in Schauder spaces, the term corresponding to the third pair of braces in the right hand side  of the equality that defines  $P_{ljr}[g,\mu]$ belongs  to $C^{m-1,\beta}(\partial\Omega)$  in case of  statement (i) and to $C^{m-1,\omega_1(\cdot)}(\partial\Omega)$ in case of  statement (ii).\par

The proof for the term corresponding to the fourth pair of braces in the right hand side  of the equality that defines  $P_{ljr}[g,\mu]$ is the same as that for the third pair.\par

By the membership of the components of $\nu$ in $C^{m,\alpha}(\partial\Omega)\subseteq
C^{m-1,\alpha}(\partial\Omega)
$, by the continuity of 
 the embedding of $C^{m,\alpha}(\partial\Omega)$ into $C^{m-1,\alpha}(\partial\Omega)$, by the continuity of the embedding of  $C^{m-1,\beta}(\partial\Omega)$ into $C^{m-2,\beta}(\partial\Omega)$,  by the continuity of
the pointwise product in Schauder spaces (cf.~\textit{e.g.}, \cite[Lems.~2.4, 2.5]{DoLa17})  and by the inductive assumption on $Q_s$, 
the term corresponding to the fifth pair of braces in the right hand side  of the equality that defines   $P_{ljr}[g,\mu]$ belongs  to $C^{m-1,\beta}(\partial\Omega)$  in case of  statement (i) and to $C^{m-1,\omega_1(\cdot)}(\partial\Omega)$ in case of  statement (ii).\par

The  membership of the components of $\nu$ in $C^{m,\alpha}(\partial\Omega)\subseteq
C^{m-1,\alpha}(\partial\Omega)
$,  the continuity of
the pointwise product in Schauder spaces (cf.~\textit{e.g.}, \cite[Lems.~2.4, 2.5]{DoLa17})  and  the continuity of the operator $v_\Omega[S_{ {\mathbf{a}} }   , \cdot]$ from the space  $C^{m-1,\beta}(\partial\Omega)$ to $C^{m,\beta}(\partial\Omega)\subseteq C^{m-1,\beta}(\partial\Omega)$   in case of  statement (i) and to $C^{m,\omega_1(\cdot)}(\partial\Omega)\subseteq C^{m-1,\omega_1(\cdot)}(\partial\Omega)$ in case of  statement (ii) 
(cf.~\cite[Thm.~7.1]{DoLa17}, Theorem \ref{thm:slay})
 imply that the term corresponding to the the last  pair of braces in the right hand side of the equality that defines  $P_{ljr}[g,\mu]$ belongs  to $C^{m-1,\beta}(\partial\Omega)$   in case of  statement (i) and to $C^{m-1,\omega_1(\cdot)}(\partial\Omega)$ in case of  statement (ii).
 Hence, the proof of (jjj) and of the theorem is complete.\hfill  $\Box$ 

\vspace{\baselineskip}

Next we  prove the following extension of a corresponding statement of \cite[Thm.~8.4]{DoLa17}.
\begin{thm}
\label{thm:rrsmn} 
Let ${\mathbf{a}}$ be as in (\ref{introd0}), (\ref{ellip}), (\ref{symr}). Let $S_{ {\mathbf{a}} }$ be a fundamental solution of $P[{\mathbf{a}},D]$. 
   Let $m\in{\mathbb{N}}\setminus\{0\}$. Let $\alpha\in]0,1]$. Let $\Omega$ be a bounded open subset of ${\mathbb{R}}^{n}$ of class $C^{m,\alpha}$.   Then the following statements hold.
\begin{enumerate}
\item[(i)] If $\alpha\in]0,1[$ and $\beta\in]0,\alpha]$, then the
 trilinear operator $R$  from the space  $\left(C^{m-1,\alpha}(\partial\Omega)\right)^{2}\times C^{m-2,\beta}(\partial\Omega)$ to $C^{m-1,\beta}(\partial\Omega)$ that is delivered by the formula
 \begin{eqnarray}\label{eq:r}
\lefteqn{
R[g,h,\mu]
\equiv \sum_{r=1}a_{r} 
\left\{
Q_r[gh,\mu]-g
Q_r[h,\mu]
-Q_r[h,g\mu]
\right\}
}
\\ \nonumber
&&\qquad\qquad\qquad 
+a\left\{
gv_\Omega [ S_{ {\mathbf{a}} }   ,h\mu  ]
- h v_\Omega [ S_{ {\mathbf{a}} }  ,g\mu  ]
\right\}\qquad\text{on}\ \partial\Omega\,
\end{eqnarray}  
 for all $(g,h,\mu)\in \left(C^{m-1,\alpha}(\partial\Omega)\right)^{2}\times C^{m-2,\beta}(\partial\Omega)$
 is continuous. 
\item[(ii)]  If $\alpha=1$ and $\beta=1$, then
the  trilinear operator $R$  from  the space $\left(C^{m-1,1}(\partial\Omega)\right)^{2}\times C^{m-2,1}(\partial\Omega)$ to $C^{m-1,\omega_1(\cdot)}(\partial\Omega)$ that is delivered by teh formula  (\ref{eq:r})  
 is continuous. 
\end{enumerate}
\end{thm}
{\bf Proof.} Since $R$ is the composition of $Q_r$ and of the single layer potential, Theorem \ref{thm:qrsmm}  on the continuity of $Q_r$ and the continuity of $v_\Omega [S_{ {\mathbf{a}} }   , \cdot]$ from $C^{m-1,\beta}(\partial\Omega)$    to $C^{m,\beta}(\partial\Omega)$   in case of  statement (i) and to $C^{m,\omega_1(\cdot)}(\partial\Omega)$ in case of  statement (ii) 
(cf.~\cite[Thm.~7.1]{DoLa17}, Theorem \ref{thm:slay})
 and the continuity of the pointwise product  
in (generalized) Schauder spaces (cf.~\textit{e.g.}, \cite[Lems.~2.4, 2.5]{DoLa17})  imply the validity of the statement.\hfill  $\Box$ 

\vspace{\baselineskip}

We are now ready to prove the following statement. 
\begin{thm}
\label{thm:wregn} 
Let ${\mathbf{a}}$ be as in (\ref{introd0}), (\ref{ellip}), (\ref{symr}). Let $S_{ {\mathbf{a}} }$ be a fundamental solution of $P[{\mathbf{a}},D]$. 
 Let $\alpha\in]0,1]$. Let $m\in{\mathbb{N}}$, $m\geq 2$. Let $\Omega$ be a bounded open subset of ${\mathbb{R}}^{n}$ of class $C^{m,\alpha}$. Then the following statements hold.
 \begin{enumerate}
\item[(i)] If $\alpha\in]0,1[$, then $W_\Omega[{\mathbf{a}},S_{{\mathbf{a}}}   ,\cdot]$
 is linear and continuous from $C^{m-1,\alpha}(\partial\Omega)$ to $C^{m,\alpha}(\partial\Omega)$.
\item[(ii)] If $\alpha=1$, then $W_\Omega[{\mathbf{a}},S_{{\mathbf{a}}}   ,\cdot]$
 is linear and continuous from $C^{m-1,1}(\partial\Omega)$ to $C^{m,\omega_1(\cdot)}(\partial\Omega)$.
 \end{enumerate}
\end{thm}
{\bf Proof.} We prove statements (i) and (ii) at the same time and make some appropriate comment when the two proofs present some difference.  By \cite[Thm.~9.1]{DoLa17}, $W_\Omega[{\mathbf{a}},S_{{\mathbf{a}}}   ,\mu]$ is continuously differentiable and the following formula holds for the tangential derivatives of $W_\Omega[{\mathbf{a}},S_{{\mathbf{a}}}   ,\mu]$ 
\begin{equation}\label{thm:wregn1}
M_{lj}[W_\Omega[{\mathbf{a}},S_{{\mathbf{a}}}   ,\mu]]=T_{lj}[\mu]
\end{equation}
where
\begin{eqnarray}\label{thm:wregn2}
\lefteqn{
T_{lj}[\mu]
\equiv W_\Omega[{\mathbf{a}},S_{{\mathbf{a}}}   ,M_{lj}[\mu]    ] 
}
\\ \nonumber
&& \qquad\quad 
+\sum_{b,r=1}^{n}a_{br}
\left\{
Q_b\left[\nu_{l},M_{jr}[\mu]\right]-
Q_b\left[\nu_{j}, M_{lr}[\mu]\right]
\right\}
\\ \nonumber
&&  \qquad\quad 
+\nu_{l} Q_j\left[ \nu\cdot a^{(1)},\mu\right] 
-\nu_{j} Q_l\left[ \nu\cdot a^{(1)},\mu\right] 
\\ \nonumber
&& \qquad\quad 
+\nu \cdot a^{(1)}
\left\{
Q_l\left[\nu_{j},\mu\right] 
-
Q_j\left[\nu_{l},\mu\right] 
\right\}
\\ \nonumber
&& \qquad\quad 
-\nu \cdot a^{(1)}
v_\Omega [S_{ {\mathbf{a}} }   , M_{lj}[\mu]]
+v_\Omega [S_{ {\mathbf{a}} }         , \nu \cdot a^{(1)}M_{lj}[\mu]]
\\ \nonumber
&&  \qquad\qquad\qquad\qquad
+R[\nu_{l},\nu_{j},\mu]
\qquad{\mathrm{on}}\ \partial\Omega\,,
\end{eqnarray}
for all $l,j\in\{1,\dots,n\}$ and $\mu\in C^1(\partial\Omega)$.

We now prove the statement by induction on $m\geq 2$. We first consider case $m=2$. 
By the definition of norm in $C^{1,\alpha}(\partial\Omega)$ and in $C^{1,\omega_1(\cdot)}(\partial\Omega)$ and by formula (\ref{thm:wregn1}) it suffices to prove that the following two statements hold.
\begin{enumerate}
\item[(j)] $W_\Omega[{\mathbf{a}},S_{{\mathbf{a}}}   ,\cdot]$ is continuous from $C^{1,\alpha}(\partial\Omega)$ to $C^{0}(\partial\Omega)$.
\item[(jj)] $T_{lj}[\cdot]$ is continuous from $C^{1,\alpha}(\partial\Omega)$ to $C^{1,\alpha}(\partial\Omega)$
 in case of  statement (i) and to $C^{1,\omega_1(\cdot)}(\partial\Omega)$ in case of  statement (ii)  for all $l$, $j\in\{1,\dots,n\}$, 
\end{enumerate}
(cf.  \cite[Lem.~2.3 (ii)]{DoLa17}). 
Since $\Omega$ is of class $C^{2,\alpha}$, then $\Omega$ is of class $C^{1,\gamma}$ for all $\gamma\in]0,1[$ and  thus  $W_\Omega[{\mathbf{a}},S_{{\mathbf{a}}}   ,\cdot]$ is continuous from $L^\infty(\partial\Omega)$ to  $C^{0}(\partial\Omega)$ (cf.~\textit{e.g.} \cite[Thm.~7.4]{DoLa17}).  Hence (j) holds true. We now prove statement (jj) by exploiting formula (\ref{thm:wregn2}). Since $\Omega$ is of class $C^{2,\alpha}$, then the normal $\nu$ belongs to 
$C^{1,\alpha}(\partial\Omega,{\mathbb{R}}^n)$. Then Theorem \ref{thm:qrsmm} with $m=2$ ensures that $Q_l\left[ \nu,\cdot\right] $ and 
$Q_j\left[ \nu\cdot a^{(1)},\cdot\right] $ are continuous from $C^{0,\alpha}(\partial\Omega)$ to $C^{1,\alpha}(\partial\Omega)$
in case of statement (i) and to $C^{1,\omega_1(\cdot)}(\partial\Omega)$ in case of statement (ii),  and that $Q_b\left[\nu_{l},M_{jr}[\cdot]\right]$ is continuous from  $C^{1,\alpha}(\partial\Omega)$ to $C^{1,\alpha}(\partial\Omega)$
in case of statement (i) and to $C^{1,\omega_1(\cdot)}(\partial\Omega)$ in case of statement (ii) for all $l$, $j$, $r\in\{1,\dots,n\}$.\par

By \cite[Thm.~1.1]{La23b}, $W_\Omega[{\mathbf{a}},S_{{\mathbf{a}}}   ,\cdot]$   is continuous from $C^{0,\alpha}(\partial\Omega)$ to $C^{1,\alpha}(\partial\Omega)$
in case of statement (i) and to $C^{1,\omega_1(\cdot)}(\partial\Omega)$ in case of statement (ii).
 Since $M_{jr}$ is continuous from $C^{1,\alpha}(\partial\Omega)$ to $C^{0,\alpha}(\partial\Omega)$, then $W_\Omega[{\mathbf{a}},S_{{\mathbf{a}}}   ,M_{lj}[\mu]    ] $ is continuous from $C^{1,\alpha}(\partial\Omega)$ to $C^{1,\alpha}(\partial\Omega)$ in case of statement (i) and to $C^{1,\omega_1(\cdot)}(\partial\Omega)$ in case of statement (ii), for all $l$, $j\in\{1,\dots,n\}$.

Since $\Omega$ is of class $C^{1,\alpha}$,  \cite[Thm.~7.1]{DoLa17} and Theorem \ref{thm:slay} imply that $V_\Omega[S_{ {\mathbf{a}} }   , \cdot]$ is continuous from $C^{0,\alpha}(\partial\Omega)$ to $C^{1,\alpha}(\partial\Omega)$
 in case of statement (i) and to $C^{1,\omega_1(\cdot)}(\partial\Omega)$ in case of statement (ii). Since $M_{lj}$ is continuous from $C^{1,\alpha}(\partial\Omega)$ to $C^{0,\alpha}(\partial\Omega)$, then $V_\Omega[S_{ {\mathbf{a}} }   , M_{lj}[\cdot]]$ is continuous from $C^{1,\alpha}(\partial\Omega)$ to $C^{1,\alpha}(\partial\Omega)$ in case of statement (i) and to $C^{1,\omega_1(\cdot)}(\partial\Omega)$ in case of statement (ii)  for all $l$, $j\in\{1,\dots,n\}$. Then the membership of $\nu$ in 
$C^{1,\alpha}(\partial\Omega,{\mathbb{R}}^n)$ and Theorem \ref{thm:rrsmn} imply that $T_{lj}$ is continuous from $C^{1,\alpha}(\partial\Omega)$ to $C^{1,\alpha}(\partial\Omega)$ in case of statement (i) and to $C^{1,\omega_1(\cdot)}(\partial\Omega)$ in case of statement (ii)  for all $l$, $j\in\{1,\dots,n\}$ and thus statement (jj) holds true. 

Hence, we have proved statements (j) and (jj) and thus 
$W_\Omega[{\mathbf{a}},S_{{\mathbf{a}}}   ,\cdot]$ is continuous from $C^{1,\alpha}(\partial\Omega)$ to $C^{2,\alpha}(\partial\Omega)$ in case of statement (i) and to $C^{2,\omega_1(\cdot)}(\partial\Omega)$ in case of statement (ii).

We now assume that $\Omega$ is of class $C^{m+1,\alpha}$ and that the statement is true for $m\geq 2$ and we turn to prove that $W_\Omega[{\mathbf{a}},S_{{\mathbf{a}}}   ,\cdot]$ is continuous from $C^{m,\alpha}(\partial\Omega)$ to $C^{m+1,\alpha}(\partial\Omega)$ in case of statement (i) and to $C^{m+1,\omega_1(\cdot)}(\partial\Omega)$ in case of statement (ii).  By the definition of norm in $C^{m+1,\alpha}(\partial\Omega)$ and in $C^{m+1,\omega_1(\cdot)}(\partial\Omega)$ and formula (\ref{thm:wregn1}), it suffices to prove that the following statements hold true.
\begin{enumerate}
\item[(a)] $W_\Omega[{\mathbf{a}},S_{{\mathbf{a}}}   ,\cdot]$ is continuous from $C^{m,\alpha}(\partial\Omega)$ to $C^{0}(\partial\Omega)$.
\item[(aa)] $T_{lj}[\cdot]$ is continuous from $C^{m,\alpha}(\partial\Omega)$ to $C^{m,\alpha}(\partial\Omega)$ in case of  statement (i) and to $C^{m,\omega_1(\cdot)}(\partial\Omega)$ in case of  statement (ii)  for all $l$, $j\in\{1,\dots,n\}$,  
\end{enumerate}
(cf. \cite[Lem.~2.3 (ii)]{DoLa17}).  
Since $C^{m,\alpha}(\partial\Omega)$ is continuously embedded into $C^{2,\alpha}(\partial\Omega)$, statement (a) follows by case $m=2$. We now prove (aa). By the inductive assumption, $W_\Omega[{\mathbf{a}},S_{{\mathbf{a}}}   ,\cdot]$ is continuous from $C^{m-1,\alpha}(\partial\Omega)$ to $C^{m,\alpha}(\partial\Omega)$ in case of  statement (i) and to $C^{m,\omega_1(\cdot)}(\partial\Omega)$ in case of  statement (ii). Since $M_{lj}[\cdot]$ is continuous from $C^{m,\alpha}(\partial\Omega)$ to $C^{m-1,\alpha}(\partial\Omega)$, we conclude that 
$W_\Omega[{\mathbf{a}},S_{{\mathbf{a}}}   ,M_{lj}[\mu]    ] $ is continuous from $C^{m,\alpha}(\partial\Omega)$ to $C^{m,\alpha}(\partial\Omega)$ in case of  statement (i) and to $C^{m,\omega_1(\cdot)}(\partial\Omega)$ in case of  statement (ii),
for all $l$, $j\in\{1,\dots,n\}$. 

Since $\Omega$ is of class $C^{m+1,\alpha}$, then the normal $\nu$ belongs to 
$C^{m,\alpha}(\partial\Omega,{\mathbb{R}}^n)$. Then Theorem \ref{thm:qrsmm}  ensures that 
$Q_l\left[ \nu,\cdot\right] $ and  $Q_r\left[ \nu\cdot a^{(1)},\cdot\right] $ are continuous from the space $C^{m-1,\alpha}(\partial\Omega)$ to $C^{m,\alpha}(\partial\Omega)$ in case of  statement (i) and to $C^{m,\omega_1(\cdot)}(\partial\Omega)$ in case of  statement (ii) and that $Q_b\left[\nu_{l},M_{jr}[\cdot]\right]$ is continuous from  $C^{m,\alpha}(\partial\Omega)$ to 
$C^{m,\alpha}(\partial\Omega)$
in case of  statement (i) and to $C^{m,\omega_1(\cdot)}(\partial\Omega)$ in case of  statement (ii)  for all $l$, $j$, $r\in\{1,\dots,n\}$.

Since $\Omega$ is of class $C^{m,\alpha}$, \cite[Thm.~7.1]{DoLa17}  and Theorem \ref{thm:slay} imply that  $V_\Omega[S_{ {\mathbf{a}} }   , \cdot]$ is continuous from $C^{m-1,\alpha}(\partial\Omega)$ to $C^{m,\alpha}(\partial\Omega)$
in case of  statement (i) and to $C^{m,\omega_1(\cdot)}(\partial\Omega)$ in case of  statement (ii). 

Since $M_{lj}$ is continuous from $C^{m,\alpha}(\partial\Omega)$ to $C^{m-1,\alpha}(\partial\Omega)$, then the operator $V_\Omega[S_{ {\mathbf{a}} }   , M_{lj}[\cdot]]$ is continuous from $C^{m,\alpha}(\partial\Omega)$ to $C^{m,\alpha}(\partial\Omega)$ in case of  statement (i) and to $C^{m,\omega_1(\cdot)}(\partial\Omega)$ in case of  statement (ii) for all $l$, $j\in\{1,\dots,n\}$. Then the membership of $\Omega$ in the class $C^{m+1,\alpha}$, of $\nu$ in 
$C^{m,\alpha}(\partial\Omega,{\mathbb{R}}^n)$  and Theorem \ref{thm:rrsmn} imply that $T_{lj}$ is continuous from $C^{m,\alpha}(\partial\Omega)$ to $C^{m,\alpha}(\partial\Omega)$ in case of  statement (i) and to $C^{m,\omega_1(\cdot)}(\partial\Omega)$ in case of  statement (ii) for all $l$, $j\in\{1,\dots,n\}$ and thus statement (aa) holds true. 

Hence, we have proved the validity of (a), (aa) and 
$W_\Omega[{\mathbf{a}},S_{{\mathbf{a}}}   ,\cdot]$ is continuous from $C^{m,\alpha}(\partial\Omega)$ to $C^{m+1,\alpha}(\partial\Omega)$ in case of  statement (i) and to $C^{m+1,\omega_1(\cdot)}(\partial\Omega)$ in case of  statement (ii) and the proof is complete. \hfill  $\Box$ 

\vspace{\baselineskip}

\section{An integral operator associated to the conormal derivative of a single layer potential}

Another relevant layer potential operator associated to the analysis of boundary value problems for the operator $P[{\mathbf{a}},D]$ is defined by
\[
W_{\ast,\Omega}[{\mathbf{a}}, S_{ {\mathbf{a}} },\mu](x)\equiv
\int_{\partial\Omega}\mu(y)DS_{ {\mathbf{a}} }(x-y)a^{(2)}\nu(x)\,d\sigma_{y}\qquad\forall x\in\partial\Omega 
\] 
for all $\mu\in C^{0}(\partial\Omega)$. We now show that Theorems   \ref{thm:slay},    \ref{thm:qrsmm}, \ref{thm:wregn}, \cite[Thm.~7.1]{DoLa17}  imply the validity of the following statement, that exploits an elementary  formula 
for $W_{\ast,\Omega}$ (cf.~\textit{e.g.}, 
 \cite[Proof of Thm.~10.1]{DoLa17}). We also mention that the following statement   extends the corresponding result of Kirsch \cite[Thm.~3.3 (b)]{Ki89} who has considered the case in which $S_{ {\mathbf{a}} }$ is the  fundamental solution of  the Helmholtz operator, $n=3$, $\alpha\in]0,1[$.
\begin{thm}
\label{v*regg}
Let ${\mathbf{a}}$ be as in (\ref{introd0}), (\ref{ellip}), (\ref{symr}). Let $S_{ {\mathbf{a}} }$ be a fundamental solution of $P[{\mathbf{a}},D]$. 
 Let $\alpha\in]0,1]$. Let $m\in{\mathbb{N}}$, $m\geq 2$. 
Let $\Omega$ be a bounded open subset of ${\mathbb{R}}^{n}$ of class $C^{m,\alpha}$. 
Then the following statements hold.
 \begin{enumerate}
\item[(i)] If $\alpha\in]0,1[$,  then the operator $W_{\ast,\Omega}[{\mathbf{a}}, S_{ {\mathbf{a}} },\cdot]$ is linear and continuous from $C^{m-2,\alpha}(\partial\Omega)$ to $C^{m-1,\alpha}(\partial\Omega)$.

\item[(ii)] If $\alpha=1$,  then the operator $W_{\ast,\Omega}[{\mathbf{a}}, S_{ {\mathbf{a}} },\cdot]$ is linear and continuous from $C^{m-2,1}(\partial\Omega)$ to $C^{m-1,\omega_1(\cdot)}(\partial\Omega)$.
 \end{enumerate}
\end{thm}
{\bf Proof.} We prove statements (i) and (ii) at the same time and make some appropriate comment when the two proofs present some difference. By a simple computation,  we have
\begin{eqnarray}
\label{v*regg1}
\lefteqn{
W_{\ast,\Omega}[{\mathbf{a}}, S_{ {\mathbf{a}} },\mu] 
=\sum_{b,r=1}^{n}a_{br}
Q_b[\nu_{r},\mu] 
}
\\ \nonumber
&&\qquad\qquad\qquad
-
W_\Omega[{\mathbf{a}}, S_{{\mathbf{a}}} ,\mu ] 
-
V_\Omega[S_{{\mathbf{a}}}       ,(a^{(1)}\nu) \mu]  
\qquad\text{on}\ \partial\Omega
\end{eqnarray}
for all  $\mu\in C^{0}(\partial\Omega)$ (cf.~\cite[(10.1)]{DoLa17}).

By the membership of the components of $\nu$ in $C^{m-1,\alpha}(\partial\Omega)$, Theorem \ref{thm:qrsmm} implies that
$Q_b[\nu_{r},\cdot]$ is continuous from $C^{m-2,\alpha}(\partial\Omega)$
to $C^{m-1,\alpha}(\partial\Omega)$ in case of  statement (i) and to $C^{m-1,\omega_1(\cdot)}(\partial\Omega)$ in case of  statement (ii) for all $r\in\{1,\dots,n\}$.

If $m=2$, then   $\Omega$ is of class $C^{2}$ and thus $W_\Omega[{\mathbf{a}}, S_{{\mathbf{a}}} ,\cdot ]$ is continuous from $C^{m-2,\alpha}(\partial\Omega)$ to $C^{m-1,\alpha}(\partial\Omega)$ in case of  statement (i) and to $C^{m-1,\omega_1(\cdot)}(\partial\Omega)$ in case of  statement (ii) (cf.~\cite[Thm.~1.1]{La23b}).

If $m>2$, then   $\Omega$ is of class $C^{m-1,\alpha}$ and $m-1\geq 2$. Then
 Theorem  \ref{thm:wregn} implies that  $W_\Omega[{\mathbf{a}}, S_{{\mathbf{a}}} ,\cdot ]$ is continuous from $C^{m-2,\alpha}(\partial\Omega)$
to $C^{m-1,\alpha}(\partial\Omega)$ in case of  statement (i) and to $C^{m-1,\omega_1(\cdot)}(\partial\Omega)$ in case of  statement (ii).

By the continuity of the pointwise product in Schauder spaces (cf.~\textit{e.g.}, \cite[Lems.~2.4, 2.5]{DoLa17}), the map from
$C^{m-2,\alpha}(\partial\Omega)$ to itself that takes $\mu$ to $(a^{(1)}\nu) \mu$
is continuous. Since $\Omega$ is of class $C^{m-1,\alpha}$, 
  \cite[Th.~7.1]{DoLa17} and Theorem \ref{thm:slay} imply that $V_\Omega[ S_{{\mathbf{a}}}      
,\cdot]$ is linear and continuous from $C^{m-2,\alpha}(\partial\Omega)$ to $C^{m-1,\alpha}(\partial\Omega)$ in case of  statement (i) and to $C^{m-1,\omega_1(\cdot)}(\partial\Omega)$ in case of  statement (ii).

Then formula (\ref{v*regg1}) implies  the validity of statement.\hfill  $\Box$ 

\vspace{\baselineskip}

\appendix

\section{Appendix: two classical lemmas}

We introduce the following two elementary technical lemmas for which we take no credit. For the convenience of the reader, we include a proof. 
\begin{lem}
\label{lem:gagrelip}
Let $n\in {\mathbb{N}}\setminus\{0,1\}$.  
Let $\Omega$ be a bounded open subset of ${\mathbb{R}}^{n}$ of class $C^1$. If $\varphi$, $\psi\in C^{0,1}(\partial\Omega)$, then
\[
\int_{\partial\Omega}M_{lr}[\varphi]\psi\,d\sigma=-
\int_{\partial\Omega}\varphi M_{lr}[\psi] \,d\sigma
\]
for all $l,r\in \{1,\dots,n\}$.
\end{lem}
{\bf Proof.}   By Mitrea,  Mitrea and Mitrea~\cite[Thm.~1.11.8]{MitMitMit22}, we have
\[
\int_{\partial\Omega}M_{lr}[\varphi\psi] \,d\sigma=0\,.
\]
 Then the statement follows by the Leibnitz rule. \hfill  $\Box$ 

\vspace{\baselineskip}

\begin{lem}\label{lem:lipgw}
 Let $n\in {\mathbb{N}}\setminus\{0,1\}$. Let $\Omega$ be a bounded open subset of ${\mathbb{R}}^{n}$ of class $C^{1}$. Let $\omega$ be a function from $[0,+\infty[$ to itself as in (\ref{om}). Then there exists $c_{\Omega,1}\in]0,+\infty[$ such that
\begin{eqnarray}\label{lem:lipgw1}
\lefteqn{
  |f(y)-f(x)-({\mathrm{grad}}_{\partial\Omega}f(x))\cdot(y-x)|
  }
  \\ \nonumber
&&\    \leq c_{\Omega,1}\left(\sup_{\partial\Omega}|f|+
\sup_{\partial\Omega}|{\mathrm{grad}}_{\partial\Omega}f|
+|{\mathrm{grad}}_{\partial\Omega}f:\partial\Omega|_{\omega(\cdot)}\right)|x-y|\omega(|x-y|)  
\end{eqnarray}
for all  $x$, $y\in\partial \Omega$ and for all $f\in C^{1,\omega(\cdot)}(\partial\Omega)$
   \textit{i.e.}, for all $f\in C^{1 }(\partial\Omega)$ such that  
   \[
|{\mathrm{grad}}_{\partial\Omega}f:\partial\Omega|_{\omega(\cdot)}
\equiv
\sup_{x,y\in \partial\Omega, x\neq y}\frac{|{\mathrm{grad}}_{\partial\Omega}f(x)-{\mathrm{grad}}_{\partial\Omega}f(y)|}{ \omega(|x-y|) }<+\infty
\,,
\]
where ${\mathrm{grad}}_{\partial\Omega}f$ denotes the tangential gradient of $f$. 
\end{lem}
{\bf Proof.}  By the Lemma  of the uniform cylinders, there exist $r_{\partial\Omega}$, $\delta\in]0,1[$ such that for each $p\in \partial\Omega$ there exist
$R_p\in O_n({\mathbb{R}})$ such that 
\[
C(p,R_p,r_{\partial\Omega},\delta)
\equiv p+ R_p^{t}({\mathbb{B}}_{n-1}(0,r_{\partial\Omega})\times]-\delta,\delta[    )
\]
is a coordinate cylinder for $\partial\Omega$ around $p$, \textit{i.e.},  
there exists a continuously differentiable function  $\gamma_p$ from ${\mathbb{B}}_{n-1}(0,r_{\partial\Omega})$ to $ ]-\delta/2 ,\delta/2 [$ such that $\gamma_p(0)=0$ and
\begin{eqnarray*}
\lefteqn{
R_p(\Omega-p )\cap ({\mathbb{B}}_{n-1}(0,r_{\partial\Omega})\times ]-\delta,\delta[)
}
\\ \nonumber
&&\qquad
=\{(\eta,y)\in {\mathbb{B}}_{n-1}(0,r_{\partial\Omega})\times]-\delta ,\delta [:\,y<\gamma_p(\eta) \}\,,
\end{eqnarray*}
and that the corresponding function $\gamma_p$  satisfies  the conditions
\[
D\gamma_p(0)=0\qquad \forall p\in \partial\Omega\,,\quad
A\equiv\sup_{p\in \partial\Omega}\|\gamma_p\|_{ C^{1}(\overline{{\mathbb{B}}_{n-1}(0,r_{\partial\Omega})}) }<+\infty 
\]
(see Dalla Riva, the author and Musolino \cite[Lem.~2.63]{DaLaMu21}).
Since $\omega$ is increasing, we have 
\[
\omega (r_{\partial\Omega}/2)\leq \omega(|x-y|)\qquad\forall  (x,y) \in  \{(\partial\Omega)^2: |x-y|\geq r_{\partial\Omega}/2\}\,.
\]
Since\begin{eqnarray*}
\lefteqn{
|f(y)-f(x)- {\mathrm{grad}}_{\partial\Omega}f(x)\cdot(y-x)|(|x-y|\omega(|x-y|))^{-1}
}
\\ \nonumber
&&\ \ 
\leq 
\left(2\sup_{\partial\Omega}|f|+\sup_{\partial\Omega}|{\mathrm{grad}}_{\partial\Omega}f| {\mathrm{diam}} (\partial\Omega)\right)
[(r_{\partial\Omega}/2)\omega(r_{\partial\Omega}/2)]^{-1} 
\end{eqnarray*}
for all $(x,y) \in  \{(\partial\Omega)^2: |x-y|\geq r_{\partial\Omega}/2\}
$, it suffices to prove (\ref{lem:lipgw1}) when 
$|x-y|<r_{\partial\Omega}/2$. Since $y\in (\partial\Omega)\cap C(x,R_x,r_{\partial\Omega},\delta)$, there exists $\eta\in {\mathbb{B}}_n(0,r)$ such that
\[
y=x+R_x^t(\eta,\gamma_x (\eta))^t\,.
\]
Then we set $\phi_{x,y}(\tau)\equiv (\tau\eta,\gamma_x(\tau\eta))$ for all $\tau\in[0,1]$. As is well known, there exists an extension $\tilde{f}\in C_c^1({\mathbb{R}}^n)$ of $f$
(cf. \textit{e.g.}, \cite[Thm.~2.85]{DaLaMu21})). 
Then we have
\begin{eqnarray*}
\lefteqn{
|f(y)-f(x)- ({\mathrm{grad}}_{\partial\Omega}f(x)) \cdot (y-x)|
}
\\ \nonumber
&& 
=|f(x+R_x^t(\eta,\gamma_x (\eta))^t)-f(x+R_x^t(0,\gamma_x (0))^t)
\\ \nonumber
&&\quad - 
{\mathrm{grad}}_{\partial\Omega}f(x+R_x^t(0,\gamma_x (0))^t)\cdot R_x^t(\eta,\gamma_x (\eta))^t|
\\ \nonumber
&& 
=\biggl|
\int_0^1 ({\mathrm{grad}}_{\partial\Omega}f(x+  R_x^t(\tau\eta,\gamma_x (\tau\eta))^t))\cdot 
(R_x^t(\eta,D\gamma_x(\tau\eta)\eta)^t)
\\ \nonumber
&&\quad 
 -  ({\mathrm{grad}}_{\partial\Omega}f(x+R_x^t(0,\gamma_x (0))^t))\cdot  (R_x^t(\eta,D\gamma_x(\tau\eta)\eta)^t)
\,d\tau
\biggr|
\\ \nonumber
&& 
\leq |{\mathrm{grad}}_{\partial\Omega}f:\partial\Omega|_{\omega(\cdot)}
\\ \nonumber
&&\qquad 
\times\int_0^1\left|
\omega(|R_x^t(\tau\eta,\gamma_x (\tau\eta))^t-
R_x^t(0,\gamma_x (0))^t
|)
\right|\,\left| \phi_{x,y}'(\tau)\right|\,d\tau  
\\ \nonumber
&& 
\leq |{\mathrm{grad}}_{\partial\Omega}f:\partial\Omega|_{\omega(\cdot)}
\sup_{\tau\in[0,1]}\omega(| (\tau\eta,\gamma_x (\tau\eta))|)
\int_0^1 \left| \phi_{x,y}'(\tau)\right|\,d\tau  
\\ \nonumber
&& 
\leq |{\mathrm{grad}}_{\partial\Omega}f:\partial\Omega|_{\omega(\cdot)}
 \omega({\mathrm{length}}(\phi_{x,y}))
{\mathrm{length}}(\phi_{x,y})
\\ \nonumber
&& 
\leq |{\mathrm{grad}}_{\partial\Omega}f:\partial\Omega|_{\omega(\cdot)}
\omega\left( |\eta|\sqrt{1+{\mathrm{Lip}}^2(\gamma_x)}\right)|\eta|\sqrt{1+{\mathrm{Lip}}^2(\gamma_x)}
\\ \nonumber
&& 
\leq |{\mathrm{grad}}_{\partial\Omega}f:\partial\Omega|_{\omega(\cdot)}
\omega\left( |x-y|\sqrt{1+{\mathrm{Lip}}^2(\gamma_x)}\right)|x-y|\sqrt{1+{\mathrm{Lip}}^2(\gamma_x)}
\\ \nonumber
&& 
\leq  |{\mathrm{grad}}_{\partial\Omega}f:\partial\Omega|_{\omega(\cdot)} (1+A^2)
\omega( |x-y|)|x-y| 
\end{eqnarray*}
(see also the last inequality of (\ref{om})).\hfill  $\Box$ 

\vspace{\baselineskip}

\noindent
{\bf Statements and Declarations}
Data sharing is not applicable to this article as no data sets were generated or analysed during the current study. This paper does not have any  conflict of interest or competing interest.

  \vspace{\baselineskip}

 \noindent
{\bf Acknowledgement}  The author  acknowledges  the support of the Research 
Project GNAMPA-INdAM   $\text{CUP}\_$E53C22001930001 `Operatori differenziali e integrali in geometria spettrale'  and is indebted to Prof. Otari Chkadua and   Prof. David Natroshvili  for a number of references and to Prof.~Paolo Luzzini and to Prof.~Paolo Musolino for a number of comments on the paper.

\noindent 
Dipartimento di Matematica `Tullio Levi-Civita',\\
Universit\`a degli Studi di Padova,\\
Via Trieste 63,\\
 I-35121 Padova,\\ 
Italy\\ 
email:mldc@math.unipd.it

\end{document}